\documentclass[a4paper, english, 11pt, twoside]{amsart}
\usepackage{babel}
\usepackage[T1]{fontenc}
\usepackage{amssymb}
\usepackage{amsmath, amsfonts, mathrsfs, amscd}
\usepackage{xcolor}
\usepackage{eucal}
\usepackage[all]{xy}
\usepackage{pifont}
\usepackage{calligra}

\newtheorem{theorem-intro}{Theorem}
\newtheorem{theorem}{Theorem}[section]
\newtheorem{lemma}[theorem]{Lemma}
\newtheorem{proposition}[theorem]{Proposition}
\newtheorem{corollary}[theorem]{Corollary}

\theoremstyle{definition}
\newtheorem{definition}[theorem]{Definition}
\newtheorem{remark}[theorem]{Remark}

\newtheorem{example}[theorem]{Example}
\newtheorem{step}{Step}
\newtheorem{paso}{Step}

\def\pf{\begin{proof}}
\def\epf{\end{proof}}

\newcommand{\Pc}{{\mathcal P}}
\newcommand{\ku}{\Bbbk}
\newcommand{\Na}{\mathbb{N}}
\newcommand{\Ent}{\mathbb{Z}}

\newcommand{\lom}{(l)!_{\omega}}
\newcommand{\nlom}{(n-l)!_{\omega}}
\newcommand{\kom}{(k)!_{\omega}}
\newcommand{\nkom}{(n-k)!_{\omega}}
\newcommand{\mh}{{\mathcal M}^H}
\newcommand{\Ss}{{\mathcal S}}

\newcommand{\qbinomk}{\binom{j}{k}_{\!\omega}}

\newcommand\Tr{\operatorname{tr}}
\newcommand\Ker{\operatorname{Ker}}
\newcommand\Ima{\operatorname{Im}}
\newcommand\id{\operatorname{id}}
\newcommand\ord{\operatorname{ord}}
\newcommand\car{\operatorname{char}}
\newcommand\Hom{\operatorname{Hom}}
\newcommand\Soc{\operatorname{Soc}}
\newcommand\Rep{\operatorname{Rep}}
\newcommand\Vect{\operatorname{Vec}}
\newcommand\Ann{\operatorname{Ann}}

\newcommand\nombre{\mathcal D}
\newcommand\hsoc{H_{\textrm{soc}}}
\newcommand\dsoc{\nombre_{\textrm{soc}}}
\newcommand\Dim{\textrm{Dim}}

\newcommand{\Oc}{{\mathcal O}}

\newcommand{\R}{{\mathcal R}}

\def\Oq{\Oc_{q}}

\renewcommand{\theequation}{\thesection.\arabic{equation}}

\newcommand{\nsubset}{\subset \hspace{-11pt} /\ }
\newcommand{\es}{\hspace{-1pt}}

\allowdisplaybreaks

\definecolor{rojo}{rgb}{1,0,0}

\DeclareMathAlphabet{\mathpzc}{OT1}{pzc}{m}{it}

\begin{document}

\title[co-Frobenius Hopf algebras]{ On two finiteness conditions for Hopf algebras with nonzero integral}

\author[Andruskiewitsch, Cuadra, Etingof]{Nicol\'as Andruskiewitsch, Juan Cuadra, Pavel Etingof}

\address{\noindent N. A.: FaMAF,
Universidad Nacional de C\'ordoba. CIEM -- CONICET. %\newline \noindent
Medina Allende s/n (5000) Ciudad Universitaria, C\'ordoba,
Argentina}
\email{andrus@famaf.unc.edu.ar}

\address{J. C.: Universidad de Almer\'\i a, Dpto. \'Algebra y An\'alisis Matem\'atico.
E04120 Almer\'\i a, Spain}
\email{jcdiaz@ual.es}

\address{P. E.: Department of Mathematics, Massachusetts Institute of Technology, Cambridge,
MA 02139, USA}
\email{etingof@math.mit.edu}

\begin{abstract}
A Hopf algebra is co-Frobenius when it has a nonzero integral. It is proved that the composition length of the indecomposable injective comodules over a co-Frobenius Hopf algebra is bounded. As a consequence, the coradical filtration of a co-Frobenius Hopf algebra is finite; this confirms a conjecture by Sorin D\u{a}sc\u{a}lescu and the first author. The proof is of categorical nature and the same result is obtained for Frobenius tensor categories of subexponential growth. A family of co-Frobenius Hopf algebras that are not of finite type over their Hopf socles is constructed, answering so in the negative another question by the same authors.
\end{abstract}

\maketitle

\vspace{-0.3cm}
\section*{Introduction}

The Haar measure on a compact group $G$ induces a linear functional $\int$ on the Hopf algebra of representative functions on $G$. The (right) invariance property of the Haar measure reads as a condition on $\int$ that can be expressed  in Hopf algebraic terms \cite[page 28]{H}. In \cite{Sw}, Sweedler extended the notion of (right invariant) integral to arbitrary Hopf algebras by means of this condition. However, not every Hopf algebra admits a nonzero (right) integral; those that do are called co-Frobenius. There is an obvious left version but a right co-Frobenius Hopf algebra is automatically left co-Frobenius.
Two main examples arose early in the study of this notion:

\begin{itemize}
  \item finite dimensional Hopf algebras \cite{LS},
  \item cosemisimple Hopf algebras \cite{Sw}.
\end{itemize}

It became slowly clear that the existence of a nonzero integral is fundamentally linked with these properties: finiteness and semisimplicity. There are two instances where these relations are apparent. The first is cohomological. Let $H$ be a Hopf algebra over a field $\ku$ and let $\mh$ denote the category of right $H$-comodules. Given $M \in \mh$, its injective hull is denoted by $E(M)$. The next characterization summarizes several results along the years, see \cite[Theorems 3 and 10]{L}, \cite[Proposition 2.3]{DN}, \cite[Lemma 1]{Don1}, \cite[page 223]{Don2} and \cite[Theorems 2.3 and 2.8]{AC}.

\begin{theorem-intro}\label{th:oldchar}
The following statements are equivalent:
\begin{enumerate} \renewcommand{\theenumi}{\roman{enumi}}\renewcommand{\labelenumi}{(\theenumi)}
\item $H$ is co-Frobenius.
\item\label{item:findim-injhull} $E(S)$ is finite dimensional for every $S \in \mh$ simple.
\item $E(\ku)$ is finite dimensional.
\item $\mh$ has a nonzero finite dimensional injective object.
\item Every $0\neq M \in \mh$ has a nonzero finite dimensional quotient.
\item $\mh$ possesses a nonzero projective object.
\item Every $M \in \mh$ has a projective cover.
\item Every injective in $\mh$ is projective.
\end{enumerate}
\end{theorem-intro}

The second instance is the heuristic principle, suggested by several examples and  results, that \emph{co-Frobenius Hopf algebras are somehow finite over a cosemisimple subobject.} The largest cosemisimple subcoalgebra of $H$ is the coradical $H_0$, the first member of the coradical filtration $(H_n)_{n\ge 0}$. The relation between the following statements was early observed by Radford in \cite[Corollary 2]{R}:

\begin{enumerate} \renewcommand{\theenumi}{\alph{enumi}}   \renewcommand{\labelenumi}{(\theenumi)}
\item $H$ is co-Frobenius.
\item The coradical filtration of $H$ is finite.
\end{enumerate}

It was proved there that (a) implies (b) under the assumption that $H_0$ is a Hopf subalgebra. This was derived from \cite[Proposition 4]{R} stating that $H=H_0E(\ku)$ for $H$ co-Frobenius. Later, it was shown in \cite[Theorem 2.1]{AD} that (b) implies (a) and
it was conjectured there that (a) $\Rightarrow$ (b) always holds. The first main result of this paper (see Section \ref{sec:finite-coradical}) is a positive answer to this conjecture, thus establishing that

\begin{theorem-intro}\label{thm:conj}
A Hopf algebra is co-Frobenius if and only if its coradical filtration is finite.
\end{theorem-intro}

The strategy of the proof is to use that the finiteness of the coradical filtration is equivalent to bound the Loewy length of all indecomposable injective objects in $\mh$. We find, more strongly, that when $H$ is co-Frobenius their composition length is bounded by $d\dim E(\ku),$ where $d$ is the largest dimension of a composition factor of $E(\ku)$, Theorem \ref{refined-bound}. The proof exploits the tensor structure of the category of finite dimensional $H$-comodules and the existence of injective hulls. Indeed, we observe in Section 2 that the same result holds for Frobenius tensor categories of subexponential growth, Theorem \ref{fi}, and provide an explicit uniform bound on the length of the indecomposable injective objects in terms of the composition series of the injective hull of the unit object. \smallbreak

In \cite{AD}, an alternative proof to Radford's theorem was given, together with an analysis of the structure of a co-Frobenius Hopf algebra whose coradical is a Hopf subalgebra, along the lines of the method proposed in \cite{AS1}; see also \cite{BDGN}.
However, there are examples of co-Frobenius Hopf algebras whose coradical is not a Hopf subalgebra.
A prominent one is the function algebra $\Oq(G)$ of a semisimple quantum group $G$ at a root of one $q$;
it was shown in \cite{APW} that the injective hulls of the simple comodules are finite dimensional.
Another approach appears in \cite{AD} through the notion of Hopf socle. Assume that $H$ has bijective antipode.
The \emph{Hopf socle} $\hsoc$ of $H$ is the span of the matrix coefficients of those simple $W \in \mh$
such that $V\otimes W$ and $W\otimes V$ are semisimple for every $V \in \mh$ simple.
If $H$ is a finitely generated module over $\hsoc$ (finite type), then $H$ is co-Frobenius \cite[Lemma 4.2]{AD}.
This is another realization of the heuristic principle above since $\Oq(G)$ is of finite type over its Hopf socle  $\Oc(G)$. The following natural question was posed in \cite[page 153]{AD}: \emph{Is any co-Frobenius Hopf algebra $H$ of finite type over $\hsoc$?} Our second main result gives a negative answer to this question. After presenting an initial direct example in Theorem \ref{motivating} and Proposition \ref{exdescrip}, we construct in Subsection \ref{subsec:blowing-dual} a new family of infinite dimensional co-Frobenius Hopf algebras $\nombre(m, \omega, (q_i)_{i\in I},\alpha)$, depending on a natural number $m$, a root of unity $\omega$ whose order $n$ divides $m$, a scalar $\alpha$, a non-empty set $I$, and a family $(q_i)_{i\in I}$ of nonzero scalars. Theorem \ref{mainth2} characterizes when $\nombre(m, \omega, (q_i)_{i\in I},\alpha)$ is of finite type over its Hopf socle and yields as a consequence:

\begin{theorem-intro}\label{question-ad}
The Hopf algebra $\nombre(m, \omega, (q_i)_{i\in I},\alpha)$ is not of finite type over its Hopf socle if $\alpha \neq 0$ and at least one of the $q_i$'s is not a root of one.
\end{theorem-intro}

The construction of $\nombre(m, \omega, (q_i)_{i\in I},\alpha)$ is inspired by the presentation by generators and
relations of the dual of a lifting of a quantum line, see Subsection \ref{subsec:dual-qline}, but blowing-up, in some sort,
part of the structure. This seems to be a novel point of view that is being explored. \smallbreak

In Section \ref{sec:positive-car} we construct other examples of co-Frobenius Hopf algebras, over fields of positive characteristic, that are not of finite type over their Hopf socles. They are smash products of a group algebra and the function algebra of a finite abelian group. An example of an infinite dimensional co-Frobenius Hopf algebra with trivial Hopf socle is given. \smallbreak

Although the answer to the question in \cite[page 153]{AD} is negative, the heuristic principle remains unscathed because the examples presented here fit into a cleft exact sequence of Hopf algebras where the kernel is finite dimensional and the cokernel cosemisimple. We wonder whether {\it any co-Frobenius Hopf algebra is an extension of some sort (short exact sequence, bosonization or else) of a finite dimensional and a cosemisimple Hopf algebra}. For example, the function algebra $\mathcal{O}(G)$ of an affine group scheme $G$ is co-Frobenius if and only if $G$ contains a linearly reductive subgroup (not necessarily normal) of finite index \cite[page 218]{Don2}. In this case, $\mathcal{O}(G)$ fits into a short exact sequence where the kernel is not a normal Hopf subalgebra but a coideal subalgebra instead.

\subsection*{Preliminaries}
For basic notions and results on Hopf algebra theory and unexplained terminology we refer to \cite{DNR,M} or \cite{Sw-book}. Throughout we will work over a ground field $\ku$. We write $\ku^{\times}$ for $\ku\backslash\{0\}.$ Vector spaces, linear maps, and unadorned tensor products are always over $\ku$. The comultiplication and counit of a coalgebra are denoted by $\Delta$ and $\varepsilon$ respectively. For a Hopf algebra $H$ its antipode is denoted by $\Ss$ and its group of group-like elements by $G(H)$. Given $g, h\in G(H)$ we set $\Pc_{g,h}(H)=\{x\in H: \Delta(x) =g \otimes x+x \otimes h\}$. A left integral $\int$ for $H$\vspace{1pt} satisfies $\int(h_{(2)})h_{(1)}= \int(h)1_H$ for all $h \in H$. Recall from \cite[page 197]{DNR} that if $\int \neq 0$ there exists a unique $g \in G(H)$ such that $\int(h_{(1)})h_{(2)}= \int(h)g$ for all $h \in H$. Such an element is called the distinguished group-like element of $H$.
\medbreak

The \emph{Loewy series} of a right $H$-comodule $M$ is the series
$$0 \subset \Soc (M) \subset \Soc^2(M) \subset \dots\subset \Soc^n(M) \subset \dots \subset \bigcup_{m \in \Na} \Soc^{m}(M)=M,$$ defined as follows: $\Soc (M)$ is the socle of $M$, i.e., the sum of all simple subcomodules of $M$. For $n>1$, $\Soc^n(M)$ is the unique subcomodule of $M$ satisfying $\Soc^{n-1}(M)\subset \Soc^{n}(M)$ and $\Soc (M/\Soc^{n-1}(M))=\Soc^n(M)/\Soc^{n-1}(M)$, see \cite[1.4]{G} or \cite[page 121]{DNR}. An alternative description of this series is through the coradical filtration: let $\rho:M \rightarrow M \otimes H$ denote the structure map of $M$, then $\Soc^{n+1}(M)=\rho^{-1}(M \otimes H_n),$ \cite[Lemma 3.1.9]{DNR}. If  $M=\Soc^n(M)$ for some $n$, the \emph{Loewy length} of $M$ is defined to be $\ell\ell(M)=\min\{m \in \Na:M=\Soc^{m}(M)\}$. Otherwise, $\ell\ell(M)=\infty$. The coradical filtration of $H$ coincides with the Loewy series of $H$, either as a right or left comodule.

\section{The coradical filtration of a co-Frobenius Hopf algebra is finite}\label{sec:finite-coradical}

 Let $H$ be a co-Frobenius Hopf algebra. To prove that $H$ has finite coradical filtration, we  use a criterium from \cite{C}.
Let $\{S_i\}_{i \in I}$ be a full set of representatives of simple right $H$-comodules. Then $H\simeq \oplus_{i\in I} E(S_i)^{n_i}$, with $n_i\in \Na$ for all $i \in I$. Since $H$ is co-Frobenius, $E(S_i)$ is finite dimensional for all $i \in I$ and hence it has finite Loewy length. Observe that $\ell\ell(E(S_i))\leq \ell(E(S_i)) \leq \dim E(S_i) $, where $\ell(E(S_i))$ denotes the composition length of $E(S_i)$. Since the Loewy series commutes with direct sums, we have:

\begin{proposition}\label{prop:cuadra-3.1}\cite[Proposition 3.1]{C} $H$ has finite coradical filtration  if and only if the set $\{\ell\ell(E(S_i))\}_{i \in I}$ is bounded. \qed
\end{proposition}

\medbreak
Thus, it would be sufficient to show that the set $\{\ell\ell(E(S_i))\}_{i \in I}$ is bounded. Indeed, we will prove the stronger statement: the set $\{\ell(E(S_i))\}_{i \in I}$ is bounded.

\begin{theorem}\label{refined-bound}
Let $S\in \mh$ be simple and $d$ the largest dimension of a composition factor of $E(\ku)$. Then $\ell(E(S)) \leq d\dim E(\ku).$
\end{theorem}

\pf Let $W$ be a composition factor of $E(S)$. Then $\Hom_{H}(E(S),E(W)) \neq 0$. Consider $E(W)$ as a subcomodule of $W \otimes E(\ku)$ and take a nonzero morphism $f: E(S) \rightarrow W \otimes E(\ku)$. We know, see for example \cite[Theorem 5.2]{C}, that $E(S)$ has a unique simple quotient, isomorphic to $\ku g \otimes S^{**}$, where $g$ is the distinguished group-like element of $H$. So $\ku g \otimes S^{**}$ is a composition factor of $\Ima(f)$ and $W \otimes E(\ku)$. There is a composition factor $U$ of $E(\ku)$ such that $\ku g \otimes S^{**}$ is a composition factor of $W \otimes U$. Then $\dim S \leq \dim W \dim U \leq d \dim W$. \par \smallbreak

Set $n= \ell(E(S))$ and let $W_1,\dots, W_n$ be the composition factors of $E(S)$. We have
\begin{align*}
\dim S \dim E(\ku) =\dim S \otimes E(\ku) \ge \dim E(S)= \sum_{j=1}^{n} \dim W_j \ge  n \frac{\dim S}{d}.
\end{align*}
From here, $n \leq d \dim E(\ku).$
\epf

A small variation of the above arguments gives a bound for the length of $S \otimes X$, with $X \in \mh$ of finite dimension, in terms of data not depending on $S$. Viewing $E(S)$ as a subcomodule of $S \otimes E(\ku)$, we get another bound for $E(S)$. It is less tight than the previous one but the proof is simpler and generalizable to tensor categories, as we will see in the next section.

\begin{proposition}\label{lema:bound-tensor}
Let $X, S\in \mh$ with $S$ simple and $\dim X < \infty$. Let $b_X$ denote the largest dimension of a composition factor of $E(\ku) \otimes X^*$.
Then $\ell(S \otimes X) \leq b_X \dim X.$
\end{proposition}

\pf Let $W$ be a composition factor of $S \otimes X$. Then $\Hom_{H}(S \otimes X, E(W)) \neq 0$. Using the adjunction, $\Hom_{H}(S,E(W) \otimes X^*) \neq 0$.
Consider $E(W)$ included in $W \otimes E(\ku)$. Since $S$ is simple, it can be viewed as a subcomodule of $W \otimes E(\ku) \otimes X^*$. There is a composition factor $U$ of $E(\ku) \otimes X^*$ such that $S$ is a composition factor of $W \otimes U$. Then $\dim S \le \dim W \dim U \leq b_X \dim W$. Now proceed as at the end of the previous proof.
\epf

This new bound can be slightly improved as follows:

\begin{corollary}\label{cor1}
Let $X, S\in \mh$ with $S$ simple and $\dim X < \infty$. Let $r$ be the number of $1$-dimensional composition factors of $X$. Then:
\begin{align*}
\ell(S \otimes X)\le b_X\dim X - r(b_X-1).
\end{align*}
\end{corollary}

\pf Set $n=\ell(X)$. Let $X_1,\ldots,X_n$ be the composition factors of $X$ and assume that $X_{n-r+1},\ldots,X_n$ are $1$-dimensional. Then $X_j^*$ is a composition factor of $X^*$ for $j=1,\ldots,n$. A composition series of $X^*$ gives rise to a series of $E(\ku) \otimes X^*$ whose factors are isomorphic to $E(\ku) \otimes X_j^*$ for $j=1,\ldots,n$. Each composition factor of $E(\ku) \otimes X_j^*$ is a composition factor of $E(\ku) \otimes X^*$. Then $b_{X_j} \le b_X$. Finally,
 \begin{align*}
\ell(S \otimes  X) & =\sum_{j=1}^n \ell(S \otimes X_j) =\sum_{j=1}^{n-r} \ell(S \otimes X_j) +  r
 \le r + \sum_{j=1}^{n-r} b_{X_j} \dim X_j \\
 & \leq b_X(\dim X-r)+r.
\end{align*}
\epf

For $X,Y \in \mh$ of finite dimension, the previous result implies that
$$\ell(Y \otimes X)\le \ell(Y)\big(b_X\dim X - r(b_X-1)\big).\vspace{5pt}$$

\begin{corollary}\label{cor:finbound}
Let $S \in \mh$ be simple. Let $b$ denote the largest dimension of a composition factor of $E(\ku) \otimes E(\ku)^*$. Assume that there are $r$ composition factors of $E(\ku)$ of dimension $1$. Then:
\begin{align}\label{eq:bound}
\ell(E(S)) \le b\dim E(\ku)-r(b-1) \leq b\dim E(\ku)-2(b-1).
\end{align}
\end{corollary}

\pf The first inequality is a consequence of Corollary \ref{cor1}. For the second one, recall from the proof of Theorem \ref{refined-bound} that $E(\ku)$ has a unique simple quotient, isomorphic to $\ku g$. Hence $r \geq 2.$ \vspace{3pt} \epf

\begin{remark}
Observe that if $H$ is pointed, then $\ell(E(S))=\ell(E(\ku))=\dim E(\ku)$ for every $S \in \mh$ simple. In this case, the bounds in Theorem \ref{refined-bound} and Corollary \ref{cor:finbound} are tight because $d=b=1$ and $r=\dim E(\ku)$. However, they are not so in general, as the next example shows.

Assume that $\ku$ contains a primitive $3$rd root of unity $q$. Denote by $H$ the dual Hopf algebra of the Frobenius-Lusztig kernel $u_q(\mathfrak{sl}_2(\ku))$. There are $3$ simple comodules: $V_0=\ku, V_1,$ and $V_2$ of dimensions 1, 2, and 3 respectively. The comodule $V_2$ is injective, see, e.g., \cite[page 158]{A}. The injective hull of $V_1$ is $V_1 \otimes V_2.$ Its composition factors are $V_0$ and $V_1$ repeated twice each. Then $\ell(E(V_1))=4$. The decomposition of $V_1 \otimes E(V_1)$ is $V_2 \oplus V_2 \oplus E(\ku)$. Hence $\dim E(\ku)=6.$ The composition factors of $E(\ku)$ are $V_0$ and $V_1$ repeated twice each. Then
$d=r=2$. The bound in Theorem \ref{refined-bound} gives $12$ to approximate $4$.

On the other hand, $E(\ku)^* \simeq E(\ku)$ because the distinguished group-like element is trivial in this example. We have the decomposition $V_1 \otimes E(\ku) \simeq V_2 \oplus V_2 \oplus E(V_1)$ and so $\ell (V_1 \otimes E(\ku))=2 +\ell(E(V_1))$. This also implies that $V_2$ occurs as a composition factor of $E(\ku) \otimes E(\ku)^*$. Hence $b=3$. The error of approximating $\ell(E(V_1))$ by $\ell (V_1 \otimes E(\ku))$ is $2$. The bound in Corollary \ref{cor:finbound} gives $14$ to approximate $4$.
\end{remark}

\section{Finiteness of the coradical filtration for Frobenius tensor categories}

We show in this section that part of the arguments used above to prove that a co-Frobenius Hopf algebra has finite coradical filtration works in the more general setting of Frobenius tensor categories of subexponential growth. This allows one to obtain the same result, for instance, for co-Frobenius co-quasi-Hopf algebras. \smallskip

We refer the reader to \cite{EO} for terminology, basic notions and the results on tensor categories needed in the sequel.

\subsection{Artinian categories and the coradical filtration}

An essentially small\footnote{This means that the isomorphism classes of objects form a set.} abelian category $\mathcal{C}$ over $\ku$ is called {\it artinian} if objects have finite length and Hom spaces are finite dimensional. This amounts to that $\mathcal{C}$ is equivalent to the category {\Large $\mathpzc{m}$}$^C$ of finite dimensional right comodules over a coalgebra $C$ over $\ku$ (which is uniquely determined up to an equivalence), see section on reconstruction theory in \cite{EGNO} and \cite[Theorem 5.1]{T}. Notice that the terminology used there for artinian is locally finite. \smallbreak

Let ${\mathcal{C}}$ be an artinian category. For $X \in {\mathcal{C}}$ denote, as before, by $\ell(X)$ and $\ell\ell(X)$ the length and Loewy length of $X$ respectively. The category ${\mathcal{C}}$ admits a filtration $({\mathcal{C}}_n)_{n \geq 0}$, called {\it coradical filtration}, where ${\mathcal{C}}_n$ is the full subcategory consisting of objects of Loewy length less or equal than $n+1$. Considering ${\mathcal{C}}$ as equivalent to {\Large $\mathpzc{m}$}$^C$, we have that ${\mathcal{C}}_n$ is equivalent to {\Large $\mathpzc{m}$}$^{C_n}$, where $C_n$ is the $n$-th member of the coradical filtration of $C$.

\subsection{Tensor categories of subexponential growth}

Assume that $\ku$ is al\-ge\-brai\-cally closed. Here, by a {\it tensor category} we mean a rigid monoidal artinian category over $\ku$ with unit object ${\mathbf 1}$, in which the tensor product is bilinear on morphisms, and ${\rm End}({\mathbf 1})=\ku$. The following definition is essentially due to Deligne, see \cite[Proposition 0.5]{De1}.

\begin{definition}
A tensor category ${\mathcal{C}}$ has subexponential growth if for any $X\in {\mathcal{C}}$ there exists a constant $K\ge 1$ such that $\ell(X^{\otimes n})\le K^n$ for large enough $n$. The infimum of such $K$ is called the spectral radius of $X$ and denoted by $\rho(X)$.
\end{definition}

\begin{remark} 1. If ${\mathcal{C}}=${\Large $\mathpzc{m}$}$^H$ for a (co-quasi-) Hopf algebra $H$, then it is clear that ${\mathcal{C}}$ has subexponential growth, and $\rho(X) \le \dim X$. \smallskip

2. More generally, assume that ${\mathcal{C}}$ admits a dimension function, i.e., a function \linebreak $X\mapsto \Dim(X) \in \Bbb R_+$ on isomorphism classes of objects, satisfying the following properties: $\Dim ({\mathbf 1})=1$, $\Dim(X^*)=\Dim(X)$, $\Dim(Z)=\Dim(X)+\Dim(Y)$ for an exact sequence $0\to X\to Z\to Y\to 0$, and $\Dim(X\otimes Y)=\Dim(X)\Dim(Y)$. Then ${\mathcal{C}}$ has subexponential growth, and $\rho(X)\le \Dim(X)$. Indeed, since for $X \ne 0$, the object $X \otimes X^*$ contains ${\mathbf 1}$, we see that $\Dim(X) \ge 1$. The additivity of $\Dim$ gives $\ell(X)\le \Dim(X)$, which implies the statement. \smallskip

3. Let $\mathcal{C}$ be a tensor category. Suppose that its Grothendieck ring admits a unital complex matrix representation $\pi$ such that $\pi(X)$ is a matrix of size $m$ with nonnegative real entries for any $X \in \mathcal{C}$ (e.g., this holds if $\mathcal{C}$ has a module category $\mathcal{M}$ with finitely many simple objects). Then $\mathcal{C}$ has subexponential growth. Indeed, we have
$$
\ell(X^{\otimes n})\le [X^{\otimes n}\otimes X^{*\otimes n}:{\mathbf 1}] \le \frac{1}{m} \Tr\big(\pi(X)^n\pi(X^*)^n\big).
$$
The latter grows exponentially with $n$ since so do the matrix elements of $\pi(X)^n$ and $\pi(X^*)^n$. Here $[X^{\otimes n}\otimes X^{*\otimes n}:{\mathbf 1}]$ denotes the multiplicity of ${\mathbf 1}$ in a composition series of $X^{\otimes n}\otimes X^{*\otimes n}$. \smallskip

4. There exist tensor categories which do not have subexponential growth, e.g., the categories ${\rm Rep}(S_t)$, $t\in \Bbb C$, obtained by
extrapolating of the representation categories of the symmetric groups $S_n$, defined by Deligne in \cite{De2}.
\end{remark}

\subsection{Frobenius tensor categories}
The following definition is inspired by the previous considerations on the category of finite dimensional comodules over a co-Frobenius Hopf algebra, see Theorem \ref{th:oldchar}.

\begin{definition}
A tensor category ${\mathcal{C}}$ is called Frobenius if each simple object has an injective hull. Equivalently, ${\mathcal{C}}$ has injective hulls.
\end{definition}

\begin{remark}
It is known, see \cite[Proposition 2.3]{EO}, that duals of projective objects in a tensor category are injective, and vice versa. Since $X \cong {}^*(X^*)$ for any object $X$, a tensor category is Frobenius if and only if it has projective covers.
\end{remark}

Semisimple tensor categories are Frobenius. In particular, Deligne's tensor category ${\rm Rep}(S_t)$ is Frobenius for $t \notin \Ent_{\ge 0}$.  Easy examples of Frobenius tensor categories that do not arise from {\Large $\mathpzc{m}$}$^H$ for a (co-quasi-) Hopf algebra $H$ can be constructed by tensoring  {\Large $\mathpzc{m}$}$^H$ with a fusion category of irrational dimension. \medskip

As before, for $X \in {\mathcal{C}}$, let $E(X)$ denote the injective hull of $X$. Let $\{S_i\}_{i \in I}$ be a full set of representatives of the simple objects in ${\mathcal{C}}$. Given $X, S \in {\mathcal{C}}$ with $S$ simple, $[X:S]$ stands for the number of occurrences of $S$ in a composition series of $X$. Notice that $[X:S]=\dim \Hom(X,E(S))$ and $[X:S] \geq \dim \Hom(S,X)$. \medskip

The next result is a generalization of Theorem \ref{refined-bound} and Corollary \ref{cor:finbound}, with a weaker bound on the length.

\begin{theorem}\label{fi}
The coradical filtration of a Frobenius tensor category ${\mathcal{C}}$ of subexponential growth is finite. More precisely, its Loewy length does not exceed
$$\sum_{i \in I}\, [E({\mathbf 1}):S_i]\, \rho(S_i\otimes E({\mathbf 1})\otimes S_i^*).$$
\end{theorem}

We need several preliminary results to establish Theorem \ref{fi}.

\begin{lemma}\label{lembd}
Let $X, S \in {\mathcal{C}}$ with $S$ simple. Then $[S \otimes X:S] \le \rho(X)$.
\end{lemma}

\pf
Write $m=[S \otimes X:S]$. Clearly, $m^r \leq [S \otimes X^{\otimes r}: S]$ for any $r \geq 1.$ Thus,
$$m^r  \le \dim \Hom(S \otimes X^{\otimes r}, E(S)) = \dim \Hom(X^{\otimes r},{}^*\hspace{-2pt}S \otimes E(S))\le K_S\, \ell(X^{\otimes r})$$
for some constant $K_S$ depending only on $S$. This implies that $m\le \rho(X)$.
\epf

\begin{proposition}\label{bd}
Let $S$ be a simple object of ${\mathcal{C}}$. Then $\ell(S \otimes X)\le \rho(X\otimes E({\mathbf 1})\otimes X^*)$.
\end{proposition}

\pf
Let $W$ be a composition factor of $S \otimes X$. Then
$$[S \otimes X: W] =\dim \Hom(S \otimes X,E(W))=\dim\Hom(S, E(W) \otimes X^*).$$
Hence $[S \otimes X: W] \leq \dim \Hom(S,W \otimes E({\mathbf 1}) \otimes X^*)$, as $E(W)$ is a direct summand of $W \otimes E({\mathbf 1})$. Let $W_1,\ldots,W_n$ be the different composition factors of $S \otimes X$. We have
$$\begin{array}{ll}
[S \otimes X \otimes E({\mathbf 1}) \otimes X^*: S] & = \displaystyle \sum_{j=1}^{n} \, [W_j \otimes E({\mathbf 1}) \otimes X^*: S][S \otimes X: W_j] \vspace{3pt} \\
 & \displaystyle \geq \sum_{j=1}^{n} \, \dim \Hom(S, W_j \otimes E({\mathbf 1}) \otimes X^*) \vspace{3pt} \\
 & \displaystyle \geq \sum_{j=1}^{n} \, \dim \Hom(S \otimes X, E(W_j)) \vspace{3pt} \\
 & \displaystyle = \sum_{j=1}^{n} \, [S \otimes X: W_j] \vspace{3pt} \\
 & = \ell (S \otimes X).
\end{array}$$
By Lemma \ref{lembd}, $\ell(S \otimes X)\le \rho(X\otimes E({\mathbf 1})\otimes X^*)$, as desired.
\epf

\begin{corollary}
For any $V\in {\mathcal{C}}$ simple,
$$\ell(E(V)) \le \sum_{i \in I}\, [E({\mathbf 1}):S_i]\,\rho(S_i\otimes E({\mathbf 1})\otimes S_i^*).$$
\end{corollary}

\pf
One has $\ell(E(V))\le \ell(V\otimes E({\mathbf 1}))= \sum_{i\in I} [E({\mathbf 1}):S_i] \ell(V\otimes S_i)$, and \vspace{1pt} the latter is bounded by
$\sum_{i\in I} [E({\mathbf 1}):S_i]\, \rho(S_i\otimes E({\mathbf 1})\otimes S_i^*)$ in view of Proposition \ref{bd}.
\epf

Now Theorem \ref{fi} follows from the fact that the length of the coradical filtration of ${\mathcal{C}}$ is the maximal Loewy length of an indecomposable  injective object.

\begin{remark}
By Proposition \ref{bd}, any Frobenius tensor category of subexponential growth satisfies $\ell(S \otimes X) \leq K(X)$, namely, $K(X)=\rho(X\otimes E({\mathbf 1})\otimes X^*)$. Conversely, any Frobenius tensor category ${\mathcal{C}}$ with this property is necessarily of subexponential growth. For, assume that for each $X,S \in {\mathcal{C}}$ with $S$ simple, there is $K(X)$ such that $\ell(S\otimes X) \le K(X)$. This implies that for $Y \in \mathcal{C}$ arbitrary, $\ell(Y \otimes X) \le \ell(Y)K(X)$. Then $\ell(X^{\otimes n})\le K(X)^n$ for all $n \geq 1$. Since Deligne's tensor category ${\rm Rep}(S_t)$ has not subexponential growth, there is no bound here for $\ell(S\otimes X)$ depending only on $X$ and not on $S$ (this is also easy to see directly, e.g., when $X$ is the analog of the permutation representation of $S_n$). \vspace{1pt}
\end{remark}

\begin{remark}
Suppose that $\mathcal{C}$ has a dimension function $\Dim$. Theorem \ref{refined-bound}, Proposition \ref{lema:bound-tensor} and Corollaries \ref{cor1} and \ref{cor:finbound} hold in $\mathcal{C}$ with exactly the same proofs. The role of the $1$-dimensional comodules is played by the invertible objects. The existence and invertibility of the distinguished group-like element $g$ is shown in \cite[2.8]{EO}. Notice that the finiteness assumption there on the isomorphism classes of simple objects is not used for this. The key point is that duals of projective objects are projective. For $S \in \mathcal{C}$ simple, that $E(S)$ has a unique simple quotient, isomorphic to $kg\otimes S^{**}$, follow from Lemmata 2.9 and 2.10 and Corollary 2.11 in \cite{EO}. With notation as there, the injective hull $E_i$ of the simple $L_i$ is isomorphic to $(P_{{}^*i})^* \simeq P_{D({}^*{i})}$. By Lemma 2.10, $P_{D(^*{i})} \simeq P_{^{**}{i}} \otimes L_g$. By Lemma 2.9 and Corollary 2.11, $P_{^{**}{i}} \simeq L_g \otimes P_{i^{**}} \otimes L_g^*$. Then $E_i \simeq L_g \otimes P_{i^{**}} \otimes L_g^* \otimes L_g \simeq L_g \otimes P_{i^{**}}$. From here, the head of $E_i$ is isomorphic to $L_g \otimes L_{i^{**}}$.
\end{remark}

\section{A family of co-Frobenius Hopf algebras not of finite type over the Hopf socle}\label{sec:answer-negative}

Before constructing the family of Hopf algebras described in the title, we present an example that tackles directly the problem of the finiteness over the Hopf socle.

\subsection{An initial example}\label{motivating}

We assume in this subsection that $\car \ku = 0$. Consider the Hopf algebra $A$ generated by $g,h$ subject to the relations $g^2=1$ and $gh=-hg$, where $g$ is  group-like and $h$ is $(g,1)$-primitive. The example will be realized as a Hopf subalgebra of the finite dual Hopf algebra $A^0$. \smallskip

For $z \in \Ent$ let $J_z$ denote the ideal generated by $h^2-z$. The algebra $A_z=A/J_z$ is isomorphic to $M_2(\ku)$ when $z \neq 0$. Writing $\bar{g},\bar{h}$ for the class of $g,h$ respectively, the isomorphism is defined by:
\begin{equation}\label{action}
\bar{g} \mapsto \left(\begin{array}{cc} 1 & 0 \\ 0 & -1 \end{array}\right), \quad
\bar{h} \mapsto \left(\begin{array}{cc} 0 & 1 \\ z & 0 \end{array}\right).
\end{equation}
Let $S_z$ be the unique (up to isomorphism) simple left $A_z$-module. It is also simple when viewed as a left $A$-module via the canonical projection $\pi_z:A \rightarrow A_z$. For $z=0$ the algebra $A_0$ is just Sweedler Hopf algebra $H_4$ and $\pi_0$ is a Hopf algebra morphism. Recall that $H_4$ has two simple modules: $\ku$, and $\ku_{\chi}$ given by the character $\chi:H_4 \rightarrow \ku, \bar{g} \mapsto -1, \bar{h} \mapsto 0$. Let $E(\ku)$ and $E(\ku_{\chi})$ denote the injective hulls (as $H_4$-modules) of $\ku$ and $\ku_{\chi}$ respectively. These four modules are (up to isomorphism) all the indecomposable left modules over $H_4$. \smallskip

Consider now the category $\Rep(A)$ of finite dimensional left $A$-modules. Let ${\mathcal C}$ be the full subcategory of $\Rep(A)$ consisting of objects $V$ on which $h^2$ acts by a semisimple linear operator with integer eigenvalues. Write $V=\oplus_{i=1}^m V(z_i)$, where $z_i \in \Ent$ and $V(z_i)=\{v \in V :h^2\cdot v=z_iv\}$ for all $i=1,...,m$. Since $h^2$ is central, $V(z_i)$ is an $A$-submodule of $V$. If $z_i \neq 0$, then the $A$-action on $V(z_i)$ factors through $A_{z_i}$ and $V(z_i) \cong S_{z_i}^{\,n_i}$ as an $A$-module for some $n_i \geq 1.$ If $z_i=0$, then the $A$-action on $V(0)$ factors through $H_4$ and $V(0)$ is isomorphic to a finite direct sum of copies of $\ku, \ku_{\chi}, E(\ku),$ and $E(\ku_{\chi})$. Conversely, $h^2$ acts by a semisimple linear operator with integer eigenvalues on any object of $\Rep(A)$ isomorphic to a finite direct sum of copies of $\ku, \ku_{\chi}, E(\ku),E(\ku_{\chi}),$ and ${S_{z_i}}'s$. This totally describes the objects of ${\mathcal C}$.

\begin{theorem}\label{initialex}
The category ${\mathcal C}$ is a tensor subcategory of $\Rep(A)$. It is tensor equivalent to {\Large $\mathpzc{m}$}$^H$ for a co-Frobenius Hopf algebra $H$ which is not of finite type over its Hopf socle.
\end{theorem}

\pf That ${\mathcal C}$ is an abelian subcategory of $\Rep(A)$ follows from the assumption on the action of $h^2$. By the same reason, $E(\ku)$ is injective in  ${\mathcal C}$, and it is the injective hull of $\ku$. Given $V, W \in {\mathcal C}$, since $h^2$ is primitive, it acts as a semisimple linear operator on $V \otimes W$. Notice that $h^2 \cdot (V(z) \otimes W(z')) \subseteq (V \otimes W)(z + z')$ for $z,z' \in \Ent$. Moreover, if $V=\oplus_{i=1}^m V(z_i)$, then  $V^*=\oplus_{i=1}^m V^*(-z_i)$ because $\Ss(h^2)=-h^2$. This shows that ${\mathcal C}$ is a tensor subcategory of $\Rep(A)$. The forgetful functor $U:{\mathcal C} \rightarrow \Vect_{\ku}$ is a fiber functor. By reconstruction theory, there is a tensor equivalence $F$ from ${\mathcal C}$ to {\Large $\mathpzc{m}$}$^H$ for some Hopf algebra $H$. This Hopf algebra must be co-Frobenius because $F(E(\ku))$ has finite length. \smallskip

Finally, we prove that $H$ is not of finite type over its Hopf socle $\hsoc$. By the form of the objects in ${\mathcal C}$, the only (up to isomorphism) simple objects are $\ku, \ku_{\chi},$ and $S_z$ for $z \in \Ent^{\diamond}$. Here $\Ent^{\diamond}=\Ent\backslash\{0\}$. It is not difficult to check directly that the multiplication rules for them are the following:
\begin{eqnarray*}
\ku_{\chi} \otimes \ku_{\chi} \simeq \ku, & \quad \ku_{\chi} \otimes S_z \simeq S_z \simeq S_z \otimes \ku_{\chi}, & \quad
S_z \otimes S_{z'} \simeq \left\{\begin{array}{ll}
H_4  & \ {\rm if}\ z' = -z, \vspace{2pt} \\
S_{z+z'}^{\,2} & \ {\rm otherwise}.
\end{array}\right.
\end{eqnarray*}
Then, $\hsoc$ has only two simple comodules and hence it is finite dimensional. Since $H$ is infinite dimensional, it cannot be of finite type over $\hsoc$. \epf

\begin{proposition}\label{exdescrip}
The Hopf algebra $H$ of the previous theorem is presented by generators $u,x,a^{\pm 1}$ and defining relations:
\begin{equation}\label{relex}
\begin{aligned}
u^2=1, & \quad x^2=0, & ux=-xu, & \quad a^{\pm 1}a^{\mp 1}=1, & ua=au, & \quad ax=xa.
\end{aligned}
\end{equation}
Its comultiplication, counit, and antipode are given by:
\begin{equation*}
\begin{aligned}
\Delta(u)& =u\otimes u, \hspace{10pt} \Delta(x)  = u \otimes x+x \otimes 1, \hspace{10pt} \Delta(a^{\pm 1})=a^{\pm 1} \otimes a^{\pm 1}\pm xua^{\pm 1} \otimes xa^{\pm 1}, \vspace{5pt} \\
\varepsilon(u)& =1, \hspace{35pt}  \varepsilon(x) =0,  \hspace{73pt} \varepsilon(a^{\pm 1}) =1, \vspace{5pt} \\
\Ss(u)& =u,  \hspace{34pt}  \Ss(x) = xu, \hspace{65pt} \Ss(a^{\pm 1})=a^{\mp 1}.
\end{aligned}
\end{equation*}
\end{proposition}

\pf Each object $X \in {\mathcal C}$ is naturally a finite dimensional right $A^0$-comodule. Let $cf(X)$ denote the coefficient space of $X$. Then $H:=\sum_{X \in {\mathcal C}} cf(X)$ is a Hopf subalgebra of $A^0$ because ${\mathcal C}$ is a tensor subcategory of $\Rep(A)$. The category {\Large $\mathpzc{m}$}$^H$ is tensor equivalent to ${\mathcal C}$. By the form of the objects in ${\mathcal C}$, it suffices to consider the family ${\mathcal F}=\{\ku, \ku_{\chi}, E(\ku), E(\ku_{\chi})\}\cup \{S_z : z \in \Ent^{\diamond}\}$ to reconstruct $H$, that is, $H=\sum_{X \in {\mathcal F}} cf(X)$. To describe the elements of $H,$ recall that $cf(X)$ is isomorphic to $(A/\Ann(X))^*$, as a coalgebra, viewing $(A/\Ann(X))^*$ inside $A^0$ through the dual map of the canonical projection of $A$ onto $A/\Ann(X)$. \smallskip

The set $\{g^ih^j: 0 \leq i \leq 1,\, 0\leq j\}$ is a basis of $A$. For each $z \in \Ent$ we consider $A_z^*$ inside $A^0$ through $\pi_z^*:A_z^* \rightarrow A^0$. Let $u,x \in H_4^*$ be defined by $\langle u, \bar{g}^i\bar{h}^j \rangle = (-1)^i\delta_{j,0}$ and $\langle x, \bar{g}^i\bar{h}^j \rangle = (-1)^i\delta_{j,1}$, with $0 \leq i,j \leq 1$. The assignment $u \rightarrow \bar{g}, x \mapsto \bar{h}$ gives a Hopf algebra isomorphism between $H_4^*$ and $H_4.$ View $u,x$ inside $A^0$ through $\pi_0^*$. Since it is a Hopf algebra morphism, we obtain the given relations and formulae for the comulplication, counit,  and antipode of $u$ and $x$. \smallskip

We now discuss the case $z \neq 0$. Recall from \eqref{action} that $A_z \simeq M_2(\ku)$. For $r,s=1,2$ let $c(z)_{rs}$ be the matrix with $1$ in the entry $(r,s)$ and zero elsewhere. Under the previous isomorphism,
\begin{equation*}
c(z)_{11} \mapsto \frac{1}{2}(\bar{1}+\bar{g}), \quad c(z)_{12} \mapsto \frac{1}{2}(\bar{h}+\bar{g}\bar{h}), \quad c(z)_{21} \mapsto \frac{1}{2z}(\bar{h}-\bar{g}\bar{h}), \quad c(z)_{22} \mapsto \frac{1}{2}(\bar{1}-\bar{g}).
\end{equation*}
Let $\{C(z)_{rs}\}_{r,s=1}^2 \subset M_2(\ku)^*$ be the dual basis of the above one. Then,
$$\Delta(C(z)_{rs})=\sum_{k=1}^2 C(z)_{rk} \otimes C(z)_{ks}, \qquad \varepsilon(C(z)_{rs})=\delta_{rs}.$$
These elements can be considered inside $A^0$ as follows:
\begin{equation*}
\begin{aligned}
\langle C(z)_{11}, g^ih^j \rangle & = z^{\frac{j}{2}} \delta_{[j],0}, & \hspace{20pt} \langle C(z)_{12}, g^ih^j \rangle & = z^{\frac{j-1}{2}} \delta_{[j],1}, \\
\langle C(z)_{21}, g^ih^j \rangle & = (-1)^iz^{\frac{j+1}{2}} \delta_{[j],1}, & \hspace{20pt} \langle C(z)_{22}, g^ih^j \rangle & = (-1)^iz^{\frac{j}{2}} \delta_{[j],0}.
\end{aligned}
\end{equation*}
Here $[j]$ stands for the class of $j$ modulo $2$. The following relations can be easily checked by direct computation:
$$\begin{array}{ll}
C(z)_{22}=uC(z)_{11}=C(z)_{11}u,   & \quad  C(z)_{21}=zxC(z)_{11}=zC(z)_{11}x, \vspace{2pt} \\
C(z)_{12}=xuC(z)_{11}=C(z)_{11}xu, & \quad C(z+z')_{11}=C(z)_{11}C(z')_{11} \ {\rm if\ } z' \neq -z, \vspace{2pt} \\
C(z)_{11}C(-z)_{11}=1, & \quad \Ss (C(z)_{11})=C(-z)_{11}.
\end{array}$$
Writing $a=C(1)_{11}$, we obtain $C(-1)_{11}=a^{-1}$ and
$$C(z)_{11}=a^{z},\ C(z)_{12}=xua^{z}, \ C(z)_{21}=zxa^{z}, \ C(z)_{22}=ua^{z}\  {\rm for\  all}\ z \in \Ent^{\diamond}.$$
Also, $ua=au$ and $xa=ax$. The comultiplication of $a$ and $a^{-1}$ read as
$$\Delta(a)=a \otimes a +xua\otimes xa, \qquad \Delta(a^{-1})=a^{-1} \otimes a^{-1}-xua^{-1} \otimes xa^{-1},$$
and the antipode $\Ss(a^{\pm 1})=a^{\mp 1}$. The previous relations show that $H$ equals the subalgebra generated by $u,x,$ and $a^{\pm 1}.$ \smallskip

We have $H=H_4+(\sum_{z \in \Ent^{\diamond}} A_z^*).$ The second sum is direct since it consists of simple subcoalgebras. Moreover,
$H_4 \cap (\sum_{z \in \Ent^{\diamond}} A_z^*)=0$ because otherwise either $\ku$ or $\ku u$ would be contained in some $A_z^*$, which is not possible. This, together with the previous relations, implies that the set $\{x^ju^ia^{z}:0\leq i,j \leq 1,\, z \in \Ent\}$ is a basis of $H$. From here, it easily follows that $H$ is presented by $u,x$ and $a^{\pm 1}$ and defining relations \eqref{relex}. Notice that $H \simeq H_4 \otimes \ku\Ent$ as algebras.
\epf \smallskip

In the next subsections we construct a large family of infinite dimensional co-Frobenius Hopf algebras that, apart from producing other examples to the above question, are interesting in its own for several reasons. This family is completely new and it does not fit in any of the general approaches to construct examples of co-Frobenius Hopf algebras, \cite[Section 3]{AD} and \cite[Section 4]{BDGN}. It provides examples of Hopf algebras generated by the coradical, a property stressed in \cite[Theorem 1.3]{AC}. We will adopt a point of view different to the above one, though related, that is susceptible of generalizations. Our construction will be better understood from the analysis of the Hopf algebra dual to a lifting of a quantum line.

\subsection{The dual of a lifting of a quantum line}\label{subsec:dual-qline}

Let $G$ be a finite abelian group and $\widehat{G}$ its group of characters. Suppose that $\car \ku \nmid \vert G \vert$. Take $1 \neq g \in G$, $\chi \in \widehat{G}$ and $\alpha \in \ku$. Let $\omega=\chi(g)$ and $n=\ord \omega$.
We assume that $n>1$ and
\begin{equation}\label{eq:hip-lift-ql}
\alpha \neq 0 \implies \chi^n=1 \text{\ and\ } g^n \neq 1.
\end{equation}
Let $G$ act on the polynomial algebra $\ku[x]$ by $\sigma\cdot x=\chi(\sigma)x$, $\sigma \in G$. The quotient of
the smash product $\ku[x]\# \ku G$ by the ideal generated by $x^n - \alpha(1-g^n)$ is a Hopf algebra, denoted by $H(G,g,\chi,\alpha)$;
it has dimension $\vert G \vert n$ and basis $\{x^j\sigma: 0\leq j <n, \sigma \in G\}$.
The elements of $G$ are group-like and $x\in \Pc_{g,1}(H(G,g,\chi,\alpha))$.
Defined in \cite[Section 5]{AS1}, it is a lifting of a quantum line in the sense of \cite[Section 4]{AS1}.

\medbreak
If $\ku$ is algebraically closed we can always take $\alpha \in \{0,1\}$ replacing $x$ by $\alpha^{- 1/n}x$ for $\alpha \neq 0$. When $G$ is cyclic of order $n$ generated by $g$ and $\alpha=0$ we get Taft Hopf algebra $T_n(\omega)$. The case $\alpha=1$ and $G$ cyclic of order $pn$ was constructed in \cite[Section 2]{Rad} to give an example of noncommutative noncocommutative Hopf algebra whose Jacobson radical is not a Hopf ideal.

\medbreak
Put $K=\Ker \chi$, $p=\vert K \vert$, and $m=\ord \chi$; $n$ divides $m$ and if $\alpha \neq 0$, then $m=n$ by assumption.
Choose $u \in G$ such that $\chi(u)=\eta$, with $\eta$ a primitive $m$-th root of unity.
The quotient group $G/K$ is cyclic of order $m$ and it is generated by the class of $u$.
We fix now a decomposition $C_{p_1} \oplus \dots  \oplus C_{p_s}$ of $K$ as a direct sum of cyclic groups of orders $p_1, \dots, p_s$.
For $i=1,\dots ,s$ let $a_i \in K$ denote a generator of the subgroup $C_{p_i}$. Every element of $G$ may be uniquely expressed as $u^ta_1^{e_1}\dots a_s^{e_s}$ with $0 \leq t <m$ and $0 \leq e_i < p_i$. We abbreviate $(e_1,\dots, e_s)$ to $e$, $a_1^{e_1}\dots a_s^{e_s}$ to $a^e$, and so on.
In particular, $g=u^{\gamma}a_1^{f_1}\dots a_s^{f_s} = u^{\gamma}a^{f}$. \emph{When $\alpha \neq 0$, we can take $u=g$ so that $\gamma=1$ and $f=0$}.  We have $\omega=\chi(g)=\chi(u^{\gamma}a^{f})=\chi(u)^{\gamma}=\eta^{\gamma}$. We also write $g^n=a^{\theta}$.
Clearly, $p=\vert K \vert=p_1\dots p_s$ and $\vert G \vert=pm$. Write $d_i=p/p_i$. Take $\xi \in \ku$ a primitive $\vert G \vert$-th root of unity such that $\xi^p=\eta$. \medbreak

Define $U,X,A_i: H(G,g,\chi,\alpha) \rightarrow \ku$ by
\begin{align*}
\langle U, x^ju^ta_1^{e_1}\dots a_s^{e_s} \rangle &= \eta^t \delta_{j,0},  \\
\langle X, x^ju^ta_1^{e_1}\dots a_s^{e_s} \rangle &= \delta_{j,1},  \\
\langle A_i, x^ju^ta_1^{e_1}\dots a_s^{e_s} \rangle &= \xi^{d_i(\theta_it+me_i)}\delta_{j,0}.
\end{align*}
Let now $D=D(G,g,\chi,\alpha)$ be the dual Hopf algebra of $H(G,g,\chi,\alpha)$.
For $0 \leq k \leq n$ recall that the $\omega$-factorial and $\omega$-binomial coefficients are given by:
\begin{align*}
(k)_{\omega} &= \sum_{j=0}^{k-1} \omega^j, & (k)!_{\omega} &= \prod_{j=1}^k (j)_{\omega},& \binom{n}{k}_{\omega} &= \frac{(n)!_{\omega}}{(k)!_{\omega} (n-k)!_{\omega}}.
\end{align*}

\begin{proposition}\label{dual}
The algebra $D$ is generated by $U$, $X$ and $A_i$, $i=1,\dots ,s$, subject to the relations:
\begin{equation}
\label{eq:relfin}
\begin{aligned}
U^m &=1, &   X^n &=0, &   A_i^{p_i}&=U^{\theta_i}, &   UX &=\omega XU,  \\
UA_i&=A_iU, &   A_iX &= \xi^{d_i(\theta_i\gamma+mf_i)}XA_i, & \ A_iA_r & =A_rA_i. &&
\end{aligned}
\end{equation}
Its comultiplication, counit, and antipode are given by $U\in G(D)$, $X\in \Pc_{U, 1}(D)$ and
\begin{align}\label{eq:delta-Ai}
\Delta(A_i) &=A_i \otimes A_i + \alpha \Big(1-\xi^{d_i\theta_im}\Big)\sum_{k=1}^{n-1} \frac{1}{\kom\nkom} X^{n-k}U^kA_i \otimes X^kA_i,
\\ \label{eq:counit-antipode-Ai}
\varepsilon(A_i)&=1, \qquad \Ss(A_i)=A_i^{p_i-1}U^{m-\theta_i}.
\end{align}
\end{proposition}
\medbreak

We stress again that $m=n$, $\gamma=1$ and $f=0$ when $\alpha \neq 0$.

\pf Throughout we shall use that $\{x^ju^ta^e\}$ is a basis of $H(G,g,\chi,\alpha)$. We divide the proof into several steps.
\begin{step}\label{step:uno}  To verify the relations in \eqref{eq:relfin} one needs the formula
\begin{align*}
\Delta(x^j\sigma)=\sum_{k=0}^j \qbinomk x^{j-k}g^k\sigma \otimes x^k \sigma.
\end{align*}
\end{step}
\noindent The computations are straightforward; along the way, one establishes the identities:
\begin{align*}
\langle U^{t'}, x^ju^ta^e \rangle &= \eta^{tt'} \delta_{j,0}; & &\text{this implies } U^m=1;
\\
\langle X^{j'}, x^ju^ta^e \rangle &= (j')!_{\omega} \delta_{j,j'}; & &\text{this implies } X^n=0;
\\
\langle A_i^{e'_i}, x^ju^ta^e \rangle &= \xi^{e'_id_i(\theta_it+me_i)}\delta_{j,0}; & &\text{this implies } A_i^{p_i}=U^{\theta_i}.
\end{align*}

\medbreak
\begin{step}\label{step:dos}
Consider the subalgebra $R$ of $D$ generated by $U$, $X$ and $A_i$, $i=1,\dots,s$. The set
\begin{align}\label{eq:basis}
\{X^{j'}U^{t'}A_1^{e'_1}\dots A_s^{e'_s}  \vert \ 0 \leq j'<n,\, 0\leq t'<m,\, 0\leq e'_i <p_i\}
\end{align}
spans $R$ and has $\vert G \vert n$ elements. We will show that it is also linearly independent and then
$R= D$. \end{step}

We claim that
\begin{equation}\label{eq:evaluation}
\langle X^{j'}U^{t'}A^{e'}, x^ju^ta^e \rangle = \bigg(\prod\limits_{i=1}^s \xi^{d_ie'_i(\theta_it+me_i)}\bigg)\eta^{tt'}(j')!_{\omega}\delta_{j',j}.
\end{equation}
For this, we easily check that $\langle A^{e'}, x^ju^ta^e \rangle=\delta_{j,0} \prod_{i=1}^s \xi^{d_ie'_i(\theta_it+me_i)}$.
Then
\begin{align*}
\langle X^{j'}U^{t'}A^{e'}, x^ju^ta^e \rangle &= \sum_{k=0}^j \qbinomk \langle X^{j'}U^{t'}, x^{j-k}
u^{\gamma k}a^{f k}u^ta^e \rangle \langle A^{e'}, x^k u^ta^e  \rangle  \\
 & \hspace{-15pt} = \bigg(\prod\limits_{i=1}^s \xi^{d_ie'_i(\theta_it+me_i)}\bigg) \langle X^{j'}U^{t'}, x^ju^ta^e \rangle  \\
 & \hspace{-15pt} = \bigg(\prod\limits_{i=1}^s \xi^{d_ie'_i(\theta_it+me_i)}\bigg)\sum_{k=0}^j \qbinomk \langle X^{j'}, x^{j-k}
u^{\gamma k}a^{f k}u^ta^e \rangle \langle U^{t'}, x^k u^ta^e  \rangle  \\
& \hspace{-15pt} = \bigg(\prod\limits_{i=1}^s \xi^{d_ie'_i(\theta_it+me_i)}\bigg)\eta^{tt'} \langle X^{j'}, x^ju^ta^e \rangle
\end{align*}
and \eqref{eq:evaluation} follows. To show linear independence of \eqref{eq:basis},
consider the following equation, where the $\lambda$'s are scalars and the limits in the sum are understood:
$$\sum_{j',t',e'} \lambda_{j',t',e'}X^{j'}U^{t'}A^{e'}=0.$$
Take $j,t$ and $e=(e_1,\dots,e_s)$ arbitraries and evaluate the previous sum at the element
$$\frac{1}{(j)!_{\omega}pm}\sum_{l=0}^{m-1} \sum_{k_1=0}^{p_1-1} \ldots \sum_{k_s=0}^{p_s-1} \eta^{-tl} \bigg(\prod_{r=1}^s \xi^{-d_re_r(mk_r+\theta_rl)}\bigg)x^ju^la_1^{k_1} \ldots a_s^{k_s};$$
we obtain
\begin{align*}
0 & = \sum_{j'} \frac{(j')!_{\omega}}{(j)!_{\omega}} \delta_{j',j} \bigg[\sum_{t'} \frac{1}{m} \sum_{l=0}^{m-1} \eta^{(t'-t)l} \bigg(\prod_{i=1}^s \xi^{-d_ie_i\theta_il}\bigg) \bigg[\sum_{e'_1,\dots,e'_s} \bigg(\prod_{i=1}^s \xi^{d_ie'_i\theta_il}\bigg) \lambda_{j',t',e'_1,\dots,e'_s} \\
 & \qquad \times \prod_{i=1}^s \bigg(\frac{1}{p_i} \sum_{k_i=0}^{p_i-1} \xi^{d_im(e'_i-e_i)k_i}\bigg)\bigg]\bigg]  = \lambda_{j,t,e}.
\end{align*}

\medbreak
\begin{step}
Let  $\R$ be the algebra presented by generators $\overline{U}$, $\overline{X}$ and $\overline{A}_i$, $i=1,\dots,s,$ with defining relations
\eqref{eq:relfin}. By Steps \ref{step:uno} and \ref{step:dos}, we have a surjective algebra morphism $\varphi: \R \to D$.  The set
$$\{\overline{X}^{j'}\overline{U}^{t'}\overline{A}_1^{e'_1}\dots \overline{A}_s^{e'_s} \ \vert \ 0 \leq j'<n,\, 0\leq t'<m,\, 0\leq e'_i <p_i\}$$
spans $\R$ and has $\vert G \vert n$ elements. Indeed, the span of this set is a left ideal of $\R$ and contains 1. Therefore, $\varphi$ is an isomorphism. \end{step}

We now proceed to establish the formulae for the comultiplication, counit, and antipode. The comultiplication at $U$ and $X$ determines the values of the counit and antipode at them, so we can skip their computations.

\medbreak
\begin{step} $\Delta(U)=U \otimes U$. \vspace{-17pt} \end{step}
$$\langle \Delta(U), (x^ju^ta^e) \otimes (x^{j'}u^{t'}a^{e'})\rangle = \langle U, (x^ju^ta^e)(x^{j'}u^{t'}a^{e'})\rangle = \eta^{tj'} \langle U, x^{j+j'}u^{t+t'}a^{e+e'}\rangle. \vspace{5pt}$$
Let $[j+j']$ be the residue class of $j+j'$ modulo $n$. We distinguish three cases:
\begin{align*}
j+j'>n &\Rightarrow & \eta^{tj'} \langle U, x^{j+j'}u^{t+t'}a^{e+e'}\rangle & = \eta^{tj'} \langle U, \alpha x^{[j+j']}(1-g^n)u^{t+t'}a^{e+e'}\rangle=0. \\[1pt]
j+j'=n &\Rightarrow & \eta^{tj'} \langle U, x^{j+j'}u^{t+t'}a^{e+e'}\rangle &  = \eta^{tj'} \langle U, \alpha (1-g^n)u^{t+t'}a^{e+e'}\rangle \\[1pt]
& &  & = \alpha\eta^{tj'} (\eta^{t+t'}-\eta^{t+t'})   =0.  \\
j+j'<n  &\Rightarrow & \eta^{tj'} \langle U, x^{j+j'}u^{t+t'}a^{e+e'}\rangle & = \eta^{tj'}\eta^{t+t'}\delta_{j+j',0}= \eta^{t+t'}\delta_{j,0}\delta_{j',0}.
\end{align*}
On the other hand,
\begin{align*}
\langle U \otimes U, (x^ju^ta^e) \otimes (x^{j'}u^{t'}a^{e'})\rangle = \langle U, x^ju^ta^e \rangle \langle U, x^{j'}u^{t'}a^{e'} \rangle
= \eta^{t+t'}\delta_{j,0}\delta_{j',0}.
\end{align*}

\bigbreak
\begin{step} $\Delta(X)=U \otimes X+X \otimes 1$.
\end{step}
Proceeding as in Step 4 one can easily check that
\begin{align*}
\langle \Delta(X), (x^ju^ta^e) \otimes (x^{j'}u^{t'}a^{e'}) \rangle & = \begin{cases} \eta^t & {\rm if\ } (j,j')=(0,1), \\
1 & {\rm if\ } (j,j')=(1,0), \\
0 & \rm{otherwise.}
\end{cases} \\
  & = \langle U \otimes X+X\otimes 1, (x^ju^ta^e) \otimes (x^{j'}u^{t'}a^{e'})\rangle.
\end{align*}

\bigbreak
\begin{step} $\Delta (A_i)$ is given by \eqref{eq:delta-Ai}.
\end{step}
We start evaluating the left-hand side of \eqref{eq:delta-Ai} at a basis element: \vspace{2pt}
\begin{align*}
\langle \Delta(A_i), x^ju^ta^e \otimes x^{j'}u^{t'}a^{e'} \rangle & = \langle A_i, (x^ju^ta^e)(x^{j'}u^{t'}a^{e'})\rangle
 = \eta^{tj'} \langle A_i, x^{j+j'}u^{t+t'}a^{e+e'}\rangle.\vspace{15pt}
\end{align*}
Again, we distinguish the three possible cases for $j+j'$:
\begin{align*}
j+j'>n &\Rightarrow & &\eta^{tj'} \langle A_i, x^{j+j'}u^{t+t'}a^{e+e'}\rangle  = \eta^{tj'} \langle A_i, \alpha x^{[j+j']}(1-g^n)u^{t+t'}a^{e+e'}\rangle   = 0. \\[1pt]
j+j'<n &\Rightarrow & &\eta^{tj'} \langle A_i, x^{j+j'}u^{t+t'}a^{e+e'}\rangle
= \eta^{tj'}\xi^{d_i(\theta_i(t+t')+m(e_i+e'_i))}\delta_{j,0}\delta_{j',0}. \\
j+j'=n &\Rightarrow & &\eta^{tj'} \langle A_i, x^{j+j'}u^{t+t'}a^{e+e'}\rangle = \eta^{tj'} \langle A_i, \alpha (1-g^n)u^{t+t'}a^{e+e'}\rangle \\
 & & & \hspace{40pt} = \eta^{tj'}\alpha \Big(\langle A_i, u^{t+t'}a^{e+e'}\rangle - \langle A_i, u^{t+t'}a^{\theta+e+e'}\rangle\Big) \\
 & & & \hspace{40pt} = \alpha\eta^{tj'}\xi^{d_i(\theta_i(t+t')+m(e_i+e'_i))}(1-\xi^{d_im\theta_i}) \\
 & & & \hspace{40pt} = \mu\eta^{tj'}\xi^{d_i(\theta_i(t+t')+m(e_i+e'_i))}.
\end{align*}
We wrote $\mu=\alpha (1-\xi^{d_i\theta_im})$ for short. We next evaluate the right-hand side of \eqref{eq:delta-Ai} at the same basis element:
\begin{align*}
& \langle A_i \otimes A_i + \mu \sum_{k=1}^{n-1} \frac{1}{\kom\nkom} X^{n-k}U^kA_i \otimes X^kA_i, (x^ju^ta^e) \otimes (x^{j'}u^{t'}a^{e'})
\rangle =  \\
& \hspace{1cm} = \langle A_i, x^ju^ta^e \rangle \langle A_i, x^{j'}u^{t'}a^{e'} \rangle \\
& \hspace{1cm} \qquad + \mu \sum_{k=1}^{n-1} \frac{1}{\kom\nkom} \langle X^{n-k}U^kA_i, x^ju^ta^e \rangle\langle X^kA_i, x^{j'}u^{t'}a^{e'} \rangle \\[3pt]
& \hspace{1cm} \overset{\eqref{eq:evaluation}} =  \xi^{d_i(\theta_it+me_i)}\xi^{d_i(\theta_it'+me'_i)}\delta_{j,0}\delta_{j',0}\\
& \hspace{1cm} \qquad +\mu \sum_{k=1}^{n-1} \frac{1}{\kom\nkom} \xi^{d_i(\theta_it+me_i)}\eta^{tk}(n-k)!_{\omega}\delta_{j,n-k}\xi^{d_i(\theta_it'+me'_i)}(k)!_{\omega}\delta_{j',k}  \\
& \hspace{1cm} = \begin{cases}
0 & {\rm if\ } j+j'>n  \\
\xi^{d_i(\theta_i(t+t')+m(e_i+e'_i))}\delta_{j,0}\delta_{j',0} & {\rm if\ } j+j'<n  \\
\mu\eta^{tj'}\xi^{d_i(\theta_i(t+t')+m(e_i+e'_i))} & {\rm if\ } j+j'=n \end{cases}
\end{align*}
Finally, we obtain the formulae \eqref{eq:counit-antipode-Ai} of the counit and antipode for $A_i$. Clearly, $\varepsilon(A_i)=1$. Since
$A_i^{\vert G \vert}=A_i^{p_id_im}=U^{\theta_id_im}=1$ we have $A_i^{\vert G \vert-1}=A_i^{p_i-1}U^{m-\theta_i}$. Then
\begin{align*}
\langle \Ss(A_i), x^ju^ta^e \rangle & =  \langle A_i, u^{-t}a^{-e} \rangle \delta_{j,0}   = \xi^{-d_i(\theta_it+me_i)} \delta_{j,0}   = \langle A_i^{\vert G \vert-1}, x^ju^ta^e \rangle.
 \end{align*}
\epf

\subsection{A new family of co-Frobenius Hopf algebras}\label{subsec:blowing-dual}

Our main construction is an infinite version of $D(G,g,\chi,\alpha)$ obtained by removing the relation $A_i^{p_i}=U^{\theta_i}$ in (\ref{eq:relfin}) and replacing $\xi^{d_i(\theta_i\gamma+mf_i)}$ in $A_iX=\xi^{d_i(\theta_i\gamma+mf_i)}XA_i$ by an arbitrary $q_i \in \ku^{\times}$, and $\xi^{d_i\theta_im}$ in \eqref{eq:delta-Ai} by $q_i^n$. Indeed, $D(G,g,\chi,\alpha)$ fits into the cleft exact sequence
$$\ku \rightarrow H(C_m,g,\chi,0) \rightarrow D(G,g,\chi,\alpha) \rightarrow \ku (C_{p_1} \oplus \dots  \oplus C_{p_s}) \rightarrow \ku.$$
The proposed changes mean replacing each $C_{p_i}$ by $\Ent$ in the cokernel and lifting the dual cocycle $\ku (C_{p_1} \oplus \dots  \oplus C_{p_s}) \rightarrow H(C_m,g,\chi,0) \otimes H(C_m,g,\chi,0)$ involved in the comultiplication. Finite dimensionality is lost but
not the co-Frobenius property because the Hopf algebra is an extension of a finite dimensional by a cosemisimple. Finiteness over the Hopf socle depends on the parameters $q_i$, as we will see in Theorem \ref{mainth2}.

\medbreak
Let $1 \neq n,m \in \Na$ be such that $n$ divides $m$. Assume that $\ku$ contains a primitive $n$-th root of unity $\omega$. Let $I$ be a non-empty set and take $q_i \in \ku^{\times}$ for each $i \in I$. Pick now $\alpha \in \ku$. Consider the $\ku$-algebra $\nombre = \nombre(m, \omega, (q_i)_{i\in I},\alpha)$ presented by generators $u,x,a_i^{\pm 1} (i \in I)$ and defining relations
\begin{align}\label{rel}
\begin{array}{llll}
u^m=1,     & \hspace{10pt}  x^n=0,         & \hspace{10pt}  a_i^{\pm 1}a_i^{\mp 1}=1,    & \hspace{10pt} ux=\omega xu, \vspace{2pt} \\
ua_i=a_iu, & \hspace{10pt} a_ix=q_ixa_i,   & \hspace{10pt}  a_ia_j=a_ja_i, & \hspace{10pt}   i,j \in I.
\end{array}
\end{align}
Fix a total order $<$ in $I$. For $r\geq 1$ set $I^{[r]} = \{(i_1,\dots ,i_r) \in I^r: i_1 < \dots < i_r\}$ and $I^{[0]}= \Ent^0 = \{0\}$. Given $F=  (i_1,\dots ,i_r)\in I^{[r]}$ and $E= (e_1,\ldots,e_r) \in \Ent^r$ we write
\begin{align}\label{eq:notation-EF}
a_F^E &= a_{i_1}^{e_1}\dots a_{i_r}^{e_r}, &
 q_F^E &=q_{i_1}^{e_1}\dots q_{i_r}^{e_r}, & a_0^0 & = 1, & q_0^0=1.
\end{align}
Thus $q_F^{nE}=q_{i_1}^{ne_1}\dots q_{i_r}^{ne_r}$. Put $\Ent^{\diamond}=\Ent \backslash \{0\}$ and $(\Ent^r)^{\diamond}=(\Ent^{\diamond})^r$. \vspace{-2pt} Let $\Gamma= \bigcup_{r\geq 1} I^{[r]} \times (\Ent^r)^{\diamond}.$ Let $C_m = \langle u\rangle$ be a cyclic  group of order $m$  and let $\chi\in \widehat{C_m}$ given by $\chi(u) = \omega$.
Then $\Ent^{(I)}$ acts on $H(C_m,u,\chi,0)$ by $a_i\cdot x = q_ix$, $a_i\cdot u = u$, $i\in I$; clearly,
$\nombre \simeq H(C_m,u,\chi,0) \# \ku \Ent^{(I)}$. Hence the set
\begin{equation}\label{basis}
B= \big\{x^su^ta_F^E : 0\leq s<n,\, 0 \leq t<m, \, (F,E) \in \Gamma \cup \{(0,0)\} \big\}
\end{equation}
is a basis of $\nombre$. Alternatively, this can be shown by applying the Diamond Lemma.

\medbreak
{\it Suppose that $m=n$ when $\alpha \neq 0$;} compare with \eqref{eq:hip-lift-ql}.
If $\ku$ is algebraically closed, for our purposes, we can take $\alpha$ in the set $\{0,1\}$.

\begin{theorem}\label{extmap} The algebra $\nombre$  bears a Hopf algebra structure uniquely defined by
\begin{equation}\label{eq:D-Hopf}
\begin{aligned}
\Delta(u)&=u\otimes u, \qquad   \Delta(x)= u \otimes x+x \otimes 1, \\
\Delta(a_i^{\pm 1})&=a_i^{\pm 1} \otimes a_i^{\pm 1}+\alpha(1-q_i^{\pm n})\sum_{k=1}^{n-1} \frac{1}{\kom\nkom} x^{n-k}u^ka_i^{\pm 1} \otimes x^ka_i^{\pm 1},\hspace{-15pt} \\
\varepsilon(u) &=1, \hspace{42pt} \varepsilon(x)=0, \hspace{57pt} \varepsilon(a_i^{\pm 1}) =1,\\
\Ss(u)& =u^{m-1}, \qquad  \Ss(x) =-u^{m-1}x, \qquad \Ss(a_i^{\pm 1})=a_i^{\mp 1}, \qquad i\in I.
\end{aligned}
\end{equation}
\end{theorem}

We split the proof into three steps.

\begin{paso}
The definitions above give rise to algebra morphisms $\Delta:\nombre \rightarrow \nombre \otimes \nombre$, $\varepsilon: \nombre \rightarrow \ku$ and to an algebra antimorphism $\Ss:\nombre \rightarrow \nombre$.
\end{paso}

\pf We must verify that these maps respect the relations (\ref{rel}) defining $\nombre$. We leave to the reader the verification for $\varepsilon$ and $\Ss$. The relations $u^m=1$, $x^n=0$, and $ux=\omega xu$ are respected by $\Delta$,  $\varepsilon$ and $\Ss$, as a particular case of the Hopf algebras defined in the previous section. The computation for the relations involving the $a_i$'s is  more involved: \vspace{-10pt}

\begin{align*}
\Delta(a_i)\Delta(a_i^{-1}) & =
\bigg(a_i \otimes a_i+\alpha(1-q_i^n)\sum_{k=1}^{n-1} \frac{1}{\kom\nkom} x^{n-k}u^ka_i \otimes x^ka_i\bigg)  \\
          &\times  \bigg(a_i^{-1} \otimes a_i^{-1}+\alpha(1-q_i^{-n})\sum_{l=1}^{n-1} \frac{1}{\lom\nlom} x^{n-l}u^la_i^{-1}\otimes x^la_i^{-1}\bigg) = \clubsuit
\end{align*}

The tensorand $x^{n-k} \otimes x^k$ in the first expression multiplies $x^{n-l} \otimes x^l$ in the second expression. This product is always zero when $k=n-1$. For $k<n-1$ the product $x^kx^l$ is nonzero only for $l=1,\dots ,n-(k+1)$. But $x^{n-k}x^{n-l}=0$ for these values of $l$. Hence the product of the two big sums is zero. Then
\begin{align*}
\clubsuit & = a_ia_i^{-1} \otimes a_ia_i^{-1} +\alpha(1-q_i^{-n})\sum_{k=1}^{n-1} \frac{1}{\kom\nkom} q_i^{n-k}x^{n-k}u^k a_i a_i^{-1} \otimes q_i^k x^k a_i a_i^{-1}  \\
 &  \qquad + \alpha(1-q_i^n)\sum_{k=1}^{n-1} \frac{1}{\kom\nkom} x^{n-k}u^ka_ia_i^{-1}\otimes x^ka_ia_i^{-1} = \Delta(1).
\end{align*}

Similarly, $\Delta(a_i^{-1})\Delta(a_i)=\Delta(1)$. Now:
\begin{align*}
\Delta(u)\Delta(a_i) &  = ua_i \otimes ua_i +\alpha(1-q_i^n)\sum_{k=1}^{n-1} \frac{1}{\kom\nkom} ux^{n-k}u^ka_i \otimes ux^ka_i  \\
& = a_iu \otimes a_iu +\alpha(1-q_i^n)\sum_{k=1}^{n-1} \frac{1}{\kom\nkom} \omega^{n-k}x^{n-k}u^ka_iu \otimes\omega^k x^ka_iu \\
& = \Delta(a_i)\Delta(u);
\end{align*}
\begin{align*}
&\Delta(a_i)\Delta(x)  =   \\
& = \bigg(a_i \otimes a_i +\alpha(1-q_i^n)\sum_{k=1}^{n-1} \frac{1}{\kom\nkom} x^{n-k}u^ka_i \otimes x^ka_i\bigg)(u \otimes x+x \otimes 1) \\
& = a_iu \otimes a_ix + a_ix \otimes a_i \\
& \qquad + \alpha(1 -q_i^n)\sum_{k=1}^{n-1} \frac{1}{\kom\nkom}(x^{n-k}u^ka_iu \otimes x^ka_ix + x^{n-k}u^ka_ix \otimes x^ka_i) \\
& = a_iu \otimes a_ix + a_ix \otimes a_i + \alpha(1-q_i^n)\sum\limits_{k=2}^{n-1} \frac{q_i}{(k-1)!_{\omega}(n-k +1)!_{\omega}} x^{n-k+1}u^{k}a_i \otimes x^{k}a_i  \\
& \qquad + \alpha(1-q_i^n)\sum_{k=2}^{n-1} \frac{q_i\omega^k}{\kom\nkom}x^{n-k+1}u^ka_i \otimes x^{k}a_i \\
& = a_iu \otimes a_ix + a_ix \otimes a_i  \\
& + \alpha(1 -q_i^n)q_i\sum_{k=2}^{n-1}\es \bigg(\es\frac{\omega^k}{\kom\nkom}\es +\es \frac{1}{(k-1)!_{\omega}(n-k+1)!_{\omega}}\es\bigg)
x^{n-k+1}u^ka_i\es\es \otimes\es x^{k}a_i\es = \spadesuit
\end{align*}
\medbreak

Observe that
\begin{align*}
\frac{\omega^k}{\kom\nkom}\es +\es \frac{1}{(k-1)!_{\omega}(n-k +1)!_{\omega}} &\es = \es \frac{1}{\kom\nkom}\es +\es \frac{\omega^{n-k + 1}} {(k-1)!_{\omega}(n-k+1)!_{\omega}}.
\end{align*}
\bigbreak

\noindent Then
\begin{align*}
 \spadesuit & = q_iua_i \otimes xa_i + q_ixa_i \otimes a_i  \\
 &  \hspace{3mm} + \alpha(1-q_i^n)q_i\hspace{-1pt}\sum_{k=2}^{n-1} \left(\hspace{-1pt}\frac{1}{\kom\nkom}\es +\es \frac{\omega^{n-k+1}}{(k-1)!_{\omega}(n-k+1)!_{\omega}}\hspace{-1pt}\right)\hspace{-3pt} x^{n-k+1}u^ka_i\es
 \otimes\es x^{k}a_i \\
 & = q_iua_i \otimes xa_i + q_ixa_i \otimes a_i + \alpha(1-q_i^n)q_i\sum_{k=2}^{n-1} \frac{1}{\kom\nkom}x^{n-k+1}u^ka_i \otimes x^{k}a_i \\
 &  \hspace{5mm} + \alpha(1-q_i^n)q_i\sum_{k=1}^{n-2} \frac{\omega^{n-k}}{\kom\nkom}x^{n-k}u^{k+1}a_i \otimes x^{k+1}a_i  \\
 & = q_iua_i \otimes xa_i + q_ixa_i \otimes a_i \\
 &  \hspace{5mm} + q_i\alpha(1-q_i^n)\sum_{k=1}^{n-1} \frac{1}{\kom\nkom} (x^{n-k+1}u^ka_i\es \otimes\es x^{k}a_i+ux^{n-k}u^{k}a_i\es \otimes\es x^{k+1}a_i) \\
 & =q_i(u \otimes x+x \otimes 1)\left(a_i \otimes a_i +\alpha(1-q_i^n)\sum_{k=1}^{n-1} \frac{1}{\kom\nkom} x^{n-k}u^{k}a_i \otimes x^{k}a_i\right)\\
 & = q_i\Delta(x)\Delta(a_i).
\end{align*}
\bigbreak

\noindent Finally,
\begin{align*}
\Delta(a_i)\Delta(a_j) & = \left(a_i \otimes a_i+\alpha(1-q_i^n)\sum_{k=1}^{n-1} \frac{1}{\kom\nkom} x^{n-k}u^ka_i \otimes x^ka_i\right)  \\
          & \qquad \times  \left(a_j \otimes a_j+\alpha(1-q_j^n)\sum_{l=1}^{n-1} \frac{1}{\lom\nlom} x^{n-l}u^la_j \otimes x^la_j \right) \\
 & \overset{\text{\ding{172}}} = a_ia_j \otimes a_ia_j +\alpha(1-q_j^n)\sum_{l=1}^{n-1} \frac{1}{\lom\nlom} q_i^{n-l}x^{n-l}u^la_ia_j \otimes q_i^lx^la_ia_j  \\
 & \qquad  + \alpha(1-q_i^n)\sum_{k=1}^{n-1} \frac{1}{\kom\nkom} x^{n-k}u^ka_ia_j\otimes x^ka_ia_j  \\
 &= a_ia_j \otimes a_ia_j + \alpha(1-(q_iq_j)^n)\sum_{l=1}^{n-1} \frac{1}{\lom\nlom} x^{n-l}u^la_ia_j \otimes x^la_ia_j
\\ &\overset{\text{\ding{173}}} = \Delta(a_j)\Delta(a_i).
\end{align*}
$\text{\ding{172}}$: The product of the two big sums is 0 as in the proof of $\Delta(a_i)\Delta(a_i^{-1})=\Delta(1)$.
\newline $\text{\ding{173}}$:  $\Delta(a_j)\Delta(a_i)$ equals the upper line with subindexes $i,j$ interchanged.
\epf \bigbreak

The second step is to give a formula for $\Delta$ evaluated at any basis element of $\nombre$, needed in Step 3 to check  the coassociativity. Recall the notation $a_F^E$ in \eqref{eq:notation-EF}.
\begin{paso}\label{comultgen}
For $s\in \Na$, $F= (i_1,\dots ,i_r) \in I^{[r]}$ and $E= (e_1,\ldots,e_r) \in \Ent^r$, we have
\begin{equation}
\begin{aligned}
\Delta(x^s&u^ta_F^E) = \sum_{l=0}^s \binom{s}{l}_{\!\omega} x^{l}u^{s-l+t}a_F^E \otimes x^{s-l}u^ta_F^E \label{comultform} \\
 & + (s)!_{\omega}\alpha(1-q_{F}^{nE}) \sum_{k=s+1}^{n-1} \frac{1}{(k)!_{\omega}(n-k+s)!_{\omega}} x^{n-k+s}u^{k+t}a_F^E  \otimes x^{k}u^ta_F^E.
\end{aligned}
\end{equation}
As a consequence, $\Delta(a_F^E)=a_F^E \otimes a_F^E$ if and only if either $\alpha=0$ or $q_F^{nE}=1.$
\end{paso}

\pf Since $\Delta$ is multiplicative and $u^t \in G(\nombre)$ commutes with $a_F^E$, it is enough to establish the formula for $x^sa_F^E$. We proceed by induction on $s$, $r$ and the exponents in $E$. Suppose that $s=0$. We leave the case $r=1$ for the reader. Let $F \in I^{[r+1]},E \in \Ent^{r+1}$. Set $F'= (i_1,\dots,i_r)$ and $E'= (e_1,\dots,e_r)$. We check the case $e_{r+1} \geq 0$, the other one being analogous.
\begin{align*}
\Delta(a_F^{E}) & = \Delta(a_{F'}^{E'})\Delta(a_{i_{r+1}}^{e_{r+1}}) \\
 & = \left( a_{F'}^{E'} \otimes a_{F'}^{E'} + \alpha (1-q_{F'}^{nE'})\sum_{k=1}^{n-1}\frac{1}{\kom\nkom} x^{n-k}u^ka_{F'}^{E'} \otimes x^{k}a_{F'}^{E'} \right) \\
 & \quad \times \left(\es a_{i_{r+1}}^{e_{r+1}} \otimes a_{i_{r+1}}^{e_{r+1}} + \alpha (1-q_{i_{r+1}}^{ne_{r+1}})\es \sum_{l=1}^{n-1}\es \frac{1}{\lom\nlom} x^{n-l}u^la_{i_{r+1}}^{e_{r+1}}\es \otimes\es x^{l}a_{i_{r+1}}^{e_{r+1}}\es \right) \\
 & \overset{\text{\ding{174}}} = a_{F}^{E} \otimes a_{F}^{E} + \alpha (1-q_{i_{r+1}}^{ne_{r+1}})\sum_{l=1}^{n-1} \frac{1}{\lom\nlom} a_{F'}^{E'}x^{n-l}u^la_{i_{r+1}}^{e_{r+1}} \otimes a_{F'}^{E'}x^{l}a_{i_{r+1}}^{e_{r+1}} \\
 & \qquad + \alpha (1-q_{F'}^{nE'})\sum_{k=1}^{n-1} \frac{1}{\kom\nkom} x^{n-k}u^k a_{F}^{E} \otimes x^{k}a_{F}^{E} \\
 & = a_{F}^{E} \otimes a_{F}^{E} + \alpha (1-q_{i_{r+1}}^{ne_{r+1}})\sum_{l=1}^{n-1} \frac{q_{F'}^{nE'}}{\lom\nlom} x^{n-l}u^l a_{F}^{E} \otimes x^{l}a_{F}^{E} \\
 &  \qquad + \alpha (1-q_{F'}^{nE'})\sum_{k=1}^{n-1} \frac{1}{\kom\nkom} x^{n-k}u^k a_{F}^{E} \otimes x^{k}a_{F}^{E} \\
 & = a_{F}^{E} \otimes a_{F}^{E} + \alpha (1-q_{F}^{nE})\sum_{k=1}^{n-1} \frac{1}{\kom\nkom} x^{n-k}u^k a_{F}^{E} \otimes x^{k}a_{F}^{E}.
\end{align*}
$\text{\ding{174}}$: The product of the two big sums is 0 as in the proof of $\Delta(a_i)\Delta(a_i^{-1})=\Delta(1)$. \smallbreak

Assume finally that the statement is proved for $s>0$. Then
\begin{multline*}
\Delta(x^{s+1}a_F^{E}) = \Delta(x)\Delta(x^sa_F^E) \\
 \hspace{1cm} = \sum_{l=0}^s \binom{s}{l}_{\!\omega} \omega^lx^{l}u^{s-l+1}a_F^E \otimes x^{s-l+1}a_F^E + \sum_{l=0}^s \binom{s}{l}_{\!\omega} x^{l+1}u^{s-l}a_F^E \otimes x^{s-l}a_F^E  \\
 \hspace{-2.5cm} + (s)!_{\omega}\alpha(1-q_{F}^{nE}) \sum_{k=s+1}^{n-1} \frac{1}{(k)!_{\omega}(n-k+s)!_{\omega}} \\
\hspace{2cm} \times\Big(\omega^{n-k+s} x^{n-k+s}u^{k+1}a_F^E \otimes x^{k+1}a_F^E + x^{n-k+s+1}u^{k}a_F^E \otimes x^{k}a_F^E\Big)  \\
\hspace{-4.7cm} = u^{s+1}a_F^E  \otimes x^{s+1}a_F^E+x^{s+1}a_F^E \otimes a_F^E \\
\hspace{-1.1cm} + \sum_{l=0}^{s-1} \bigg(\omega^{l+1} \binom{s}{l+1}_{\!\omega}  + \binom{s}{l}_{\!\omega} \bigg) x^{l+1}u^{s-l}a_F^E \otimes x^{s-l}a_F^E \\
\hspace{1.8cm} + (s)!_{\omega}\alpha(1-q_{F}^{nE})\bigg(\es \sum_{k=s+1}^{n-2} \frac{\omega^{n-k+s}}{(k)!_{\omega}(n-k+s)!_{\omega}} x^{n-k+s}u^{k+1}a_F^E\es  \otimes\es x^{k+1}a_F^E\es\bigg) \\
\hspace{1.6cm} + (s)!_{\omega}\alpha(1-q_{F}^{nE})\bigg(\sum_{k=s+2}^{n-1} \frac{1}{(k)!_{\omega}(n-k+s)!_{\omega}}x^{n-k+s+1}u^{k}a_F^E  \otimes x^{k}a_F^E\bigg) \\ \hspace{-1.1cm}= \sum_{l=0}^{s+1} \binom{s+1}{l}_{\!\omega} x^{l}u^{s+1-l}a_F^E \otimes x^{s+1-l}a_F^E
 + (s)!_{\omega}\alpha(1-q_{F}^{nE})\\
\hspace{0.2cm} \times\es\es \sum_{k=s+2}^{n-1}\es\es\es \bigg(\es \frac{\omega^{n-k+s+1}}{(k-1)!_{\omega}(n-k+s+1)!_{\omega}}\es +\es   \frac{1}{(k)!_{\omega}(n-k+s)!_{\omega}}\es\bigg)\es x^{n-k+s+1}\es u^{k}a_F^E \otimes x^{k}a_F^E \\
\hspace{-4cm} \overset{\text{\ding{175}}} = \sum_{l=0}^{s+1} \binom{s+1}{l}_{\!\omega} x^{l}u^{s+1-l}a_F^E \otimes x^{s+1-l}a_F^E  \\
\hspace{1.8cm} + (s+1)!_{\omega}\alpha(1-q_{F}^{nE})\sum_{k=s+2}^{n-1} \frac{1}{\kom (n-k+s+1)!_{\omega}} x^{n-k+s+1}u^{k}a_F^E\otimes x^{k}a_F^E.
\end{multline*}

\medbreak
$\text{\ding{175}}$: We have

\medbreak
$\displaystyle\frac{\omega^{n-k+s+1}}{(k-1)!_{\omega}(n-k+s+1)!_{\omega}}+\frac{1}{(k)!_{\omega}(n-k+s)!_{\omega}}=\frac{(s+1)_{\omega}}{\kom (n-k+s+1)!_{\omega}}$.
\epf \bigbreak

\begin{paso}
The maps $\Delta$, $\varepsilon$ and $\Ss$ defined in \eqref{eq:D-Hopf} equip $\nombre$ with a Hopf algebra structure.
\end{paso}

\pf We first prove that $(\nombre,\Delta,\varepsilon)$ is a coalgebra. Since $\Delta$ and $\varepsilon$ are algebra morphisms, it suffices to check the coassociativity and counit axioms for the generators $u,x,a_i^{\pm 1}$. Clearly, it holds for $u,x$ because the algebra they generate is a particular case of the Hopf algebras discussed in the previous subsection. So, we only must check them for $a_i^{\pm 1}$. For $\alpha=0$ the verification is straightforward. We assume that $\alpha \neq 0$. Using (\ref{comultform}), we compute:

\renewcommand{\theequation}{\arabic{equation}}
\setcounter{equation}{0}

\vspace{-0.7cm}
\begin{align}
\notag
(\Delta \otimes \id) \Delta(a_i) = \Delta(a_i)\otimes a_i+\alpha(1-q_i^{n})\sum_{k=1}^{n-1} \frac{1}{\kom\nkom} \Delta(x^{n-k}u^ka_i) \otimes x^k a_i  \hspace{0.4cm} \\
\nonumber
= a_i \otimes a_i \otimes a_i + \alpha(1-q_i^{n})\sum_{l=1}^{n-1} \frac{1}{\lom\nlom} x^{n-l}u^la_i \otimes x^{l}a_i \otimes a_i
\hspace{1.9cm}
\\ \label{eq:2.5}
+ \alpha(1-q_i^{n})\sum_{k=1}^{n-1} \sum_{l=0}^{n-k} \frac{1}{\kom\nkom} \binom{n-k}{l}_{\!\! \omega}  x^lu^{n-l}a_i \otimes
x^{n-k-l}u^ka_i \otimes x^{k}a_i \\[-5pt]
+ (\alpha(1-q_i^{n}))^{2}\hspace{-0.5pt} \sum_{k=1}^{n-1} \sum_{v=n-k+1}^{n-1}\hspace{-1pt}\frac{1}{\kom (v)!_{\omega}(2n-v-k)!_{\omega}} x^{2n-v-k}u^{v+k}a_i \hspace{-1pt}\label{eq:2.6}
\otimes \hspace{-1pt} x^{v}u^ka_i \hspace{-1pt} \otimes \hspace{-1pt} x^k a_i
\end{align}

\noindent Consider the sum in \eqref{eq:2.5}. Take out the part corresponding to $l=0$. Observe now that $l$ takes all values from $1$ to $n-1$ when $k$ runs. For $l=j$ the tensorand accompanying $x^ju^{n-j}a_i$ in \eqref{eq:2.5} is:
$$\alpha(1-q_i^{n})\sum_{k=1}^{n-j} \frac{1}{\kom\nkom}\binom{n-k}{j}_{\!\! \omega} x^{n-k-j}u^ka_i \otimes x^{k}a_i.$$
Set $l=v+k-n$ in \eqref{eq:2.6}. {\it Notice that $u^n=1$ and $u^{v+k}=u^{n+l-k+k}=u^l$ because we are assuming $m=n$ for $\alpha \neq 0$}.
We continue our computation by making these substitutions: \vspace{-0.6cm}

\begin{align} \notag
= a_i \otimes a_i \otimes a_i + \alpha(1-q_i^{n})\sum_{k=1}^{n-1} \frac{1}{\kom\nkom} x^{n-k}u^ka_i \otimes x^{k}a_i \otimes a_i \hspace{1.6cm}\\ \notag + \alpha(1-q_i^{n})\sum_{k=1}^{n-1} \frac{1}{\kom\nkom} a_i \otimes x^{n-k}u^ka_i \otimes x^ka_i \hspace{3.5cm}  \\ \label{eq:2.7}
\hspace{0.7cm} + \alpha(1-q_i^{n})\sum_{l=1}^{n-1} \sum_{k=1}^{n-l} \frac{1}{\kom\nkom}\binom{n-k}{l}_{\!\! \omega}  x^lu^{n-l}a_i \otimes
x^{n-k-l}u^ka_i \otimes x^{k}a_i \hspace{-3pt} \\ \label{eq:2.8}
+ (\alpha(1\hspace{-0.7pt}-\hspace{-0.5pt}q_i^{n}))^{2}\hspace{-2pt}\sum_{k=1}^{n-1}\sum_{l=1}^{k-1} \frac{1}{\kom (n-k+l)!_{\omega}\nlom} x^{n-l}u^{l}a_i\hspace{-2pt} \otimes\hspace{-1pt} x^{n\hspace{-0.8pt}-\hspace{-0.8pt}k\hspace{-0.8pt}+\hspace{-0.8pt}l}u^ka_i \hspace{-2pt} \otimes \hspace{-1pt} x^k a_i. \hspace{-15pt}
\end{align}
In \eqref{eq:2.7} put $t=n-l$. In \eqref{eq:2.8} observe that $l$ takes all values from $1$ to $n-2$ when $k$ runs and for $l=j$ the tensorand accompanying $x^{n-j}u^{j}a_i$ is:
$$\sum_{k=j+1}^{n-2} \frac{1}{\kom (n-k+j)!_{\omega}(n-j)!_{\omega}} x^{n-k+j}u^ka_i \otimes x^k a_i.$$
Making these two substitutions we have: \vspace{-10pt}

\begin{align}\notag
= a_i \otimes a_i \otimes a_i + \alpha(1-q_i^{n})\sum_{k=1}^{n-1} \frac{1}{\kom\nkom} x^{n-k}u^ka_i \otimes x^{k}a_i \otimes a_i \hspace{1.7cm} \\ \notag
+ \alpha(1-q_i^{n})\sum_{k=1}^{n-1} \frac{1}{\kom\nkom} a_i \otimes x^{n-k}u^ka_i \otimes x^ka_i \hspace{3.7cm} \\ \label{eq:2.9}
+ \alpha(1-q_i^{n})\sum_{t=1}^{n-1} \sum_{k=1}^{t} \frac{1}{\kom\nkom} \binom{n-k}{n-t}_{\!\! \omega}  x^{n-t}u^{t}a_i \otimes x^{t-k}u^ka_i \otimes x^{k}a_i \hspace{0.4cm} \\ \nonumber
+ (\alpha(1-q_i^{n}))^{2}\es\es \sum_{l=1}^{n-1}\sum_{k=l+1}^{n-1}\es \frac{\lom}{\kom (n-k+l)!_{\omega}\nlom\lom}\es x^{n-l}u^{l}a_i\es\es \otimes\es x^{n-k+l}u^ka_i\es\es \otimes\es x^k \es a_i.\hspace{-6pt} \vspace{12pt}
\end{align}
In \eqref{eq:2.9} put $l=t-k$. We obtain: \vspace{-10pt}

\begin{align}\notag
= a_i \otimes a_i \otimes a_i+ \alpha(1-q_i^{n})\sum_{k=1}^{n-1} \frac{1}{\kom\nkom} a_i \otimes x^{n-k}u^ka_i \otimes x^ka_i \hspace{1.7cm} \\ \label{eq:2.11}
+ \alpha(1-q_i^{n})\sum_{k=1}^{n-1} \frac{1}{\kom\nkom} x^{n-k}u^ka_i \otimes x^{k}a_i \otimes a_i \hspace{3.7cm}\\ \label{eq:2.12}
+ \alpha(1-q_i^{n})\es\es \sum_{t=1}^{n-1} \sum_{l=0}^{t-1} \frac{1}{(t-l)!_{\omega}(n-t+l)!_{\omega}}\binom{n-t+l}{n-t}_{\!\! \omega}\es  x^{n-t}u^{t}a_i\es\es \otimes\es x^{l}u^{t-l}a_i\es\es \otimes\es x^{t-l}a_i \hspace{-2pt} \\ \notag
+ (\alpha(1-q_i^{n}))^{2}\es\es \sum_{k=1}^{n-1}\sum_{l=k+1}^{n-1}\es \frac{\kom}{\lom (n-l+k)!_{\omega}\nkom\kom} x^{n-k}u^{k}a_i\es\es \otimes\es x^{n-l+k}u^la_i\es\es \otimes \es x^l a_i\hspace{-2pt} \vspace{12pt}
\end{align}
We rewrite the coefficient in \eqref{eq:2.12}:
\begin{align*}
\frac{1}{(t-l)!_{\omega}(n-t+l)!_{\omega}} \binom{n-t+l}{n-t}_{\!\! \omega} & = \frac{1}{t!_{\omega}(n-t)!_{\omega}} \binom{t}{l}_{\!\! \omega}.
\end{align*}
Observe that the formula obtained in \eqref{eq:2.13} gives \eqref{eq:2.11} for $l=k$:
\begin{align}
\notag
= a_i \otimes a_i \otimes a_i + \alpha(1-q_i^{n})\sum_{k=1}^{n-1} \frac{1}{\kom\nkom} a_i \otimes x^{n-k}u^ka_i \otimes x^ka_i
\hspace{1.3cm} \\ \label{eq:2.13}
+ \alpha(1-q_i^{n})\sum_{k=1}^{n-1} \sum_{l=0}^{k} \frac{1}{(k)!_{\omega}(n-k)!_{\omega}} \binom{k}{l}_{\!\! \omega}  x^{n-k}u^{k}a_i \otimes x^{l}u^{k-l}a_i \otimes x^{k-l}a_i \hspace{0.4cm} \\ \notag
+ (\alpha(1-q_i^{n}))^{2}\es\es \sum_{k=1}^{n-1}\sum_{l=k+1}^{n-1}\es \frac{\kom}{\kom\nkom\lom (n-l+k)!_{\omega}}\es x^{n-k}u^{k}a_i\es\es \otimes \es x^{n-l+k}u^la_i\es\es \otimes\es x^l a_i \hspace{-10pt} \\ \notag
= a_i \otimes \Delta(a_i) + \alpha(1-q_i^{n})\sum_{k=1}^{n-1} \frac{1}{\kom\nkom} x^{n-k}u^ka_i \otimes \Delta(x^{k}a_i) \hspace{1.6cm} \\ \nonumber
= (\id \otimes \Delta)\Delta(a_i). \hspace{9.5cm}
\end{align}
It follows at once that $(\Delta \otimes \id)\Delta(a_i^{-1})=(\id \otimes \Delta)\Delta(a_i^{-1})$.
We leave to the reader to check the counit axiom.
We finally prove that $\Ss$ is the inverse of $\id$ for the convolution product;
it is enough to check the axioms for the generators  $a_i^{\pm 1}$. We compute:
\begin{align*}
(\Ss * \id)(a_i) & = \Ss(a_i)a_i +\alpha(1-q_i^{n})\sum\limits_{k=1}^{n-1} \frac{1}{\kom\nkom} \Ss(x^{n-k}u^ka_i)x^{k}a_i  \\
& = a_i^{-1}a_i\es + \es \alpha(1-q_i^{n})\es\es\sum_{k=1}^{n-1} \frac{(-1)^{n-k}\omega^{-\frac{(n-k-1)(n-k)}{2}}}{\kom\nkom} a_i^{-1}u^{n-k}u^{n-n+k}x^{n-k}x^{k}a_i  \\
& = 1 + \alpha(1-q_i^{n})\sum_{k=1}^{n-1} \frac{(-1)^{n-k}\omega^{-\frac{(n-k-1)(n-k)}{2}}}{\kom\nkom} a_i^{-1}x^na_i   = \varepsilon(a_i). \\
(id * \Ss)(a_i) & = a_i\Ss(a_i) + \alpha(1-q_i^n)\sum\limits_{k=1}^{n-1} \frac{1}{\kom\nkom} x^{n-k}u^ka_i\Ss(x^ka_i)  \\
& = a_ia_i^{-1} + \alpha(1-q_i^{n}) \sum\limits_{k=1}^{n-1} \frac{(-1)^k \omega^{-\frac{(k-1)k}{2}}}{\kom\nkom} x^{n-k}u^ka_ia_i^{-1}u^{n-k}x^k\\
& = 1 + \alpha(1-q_i^{n}) \sum\limits_{k=1}^{n-1} \frac{(-1)^k \omega^{-\frac{(k-1)k}{2}}}{\kom\nkom} x^n   = \varepsilon(a_i).
\end{align*}
The computation for $a_i^{-1}$ is the same replacing $a_i$ by $a_i^{-1}$ and $q_i^n$ by $q_i^{-n}$.\epf

\begin{remark}
Consider the example in the above family corresponding to $n=m=2$ and $\vert I \vert=1$. Write $q$ instead of $q_i$. If we take $\alpha=(2(1-q))^{-1}$ and make $q$ to tend to $1$, we recover the example of Subsection \ref{motivating}. A similar example can be constructed by taking $n=m, q_i=q$ for all $i \in I, \alpha=(n(1-q))^{-1}$ and making $q$ to tend to $1$. The defining relations, comultiplication, counit, and antipode in this case read as:
\begin{equation*}
\begin{aligned}
u^n &=1, & \hspace{10pt}  x^n &=0,&    \hspace{10pt}  ux  &=\omega xu,  & \hspace{10pt} a_i^{\pm 1}a_i^{\mp 1}&=1,  \\
ua_i &=a_iu, & \hspace{10pt} a_ix &=xa_i,  & \hspace{10pt} a_ia_j &=a_ja_i, & \hspace{10pt} & \vspace{7pt}
\end{aligned}
\end{equation*}
\begin{equation*}
\begin{aligned}
\Delta(u)=u\otimes u, \qquad   \Delta(x)= u \otimes x+x \otimes 1, \hspace{2.8cm} \\
\Delta(a_i^{\pm 1})=a_i^{\pm 1} \otimes a_i^{\pm 1}\pm \sum_{k=1}^{n-1} \frac{1}{\kom\nkom} x^{n-k}u^ka_i^{\pm 1} \otimes x^ka_i^{\pm 1}, \hspace{-0.6cm} \\
\varepsilon(u)=1, \hspace{42pt} \varepsilon(x)=0, \hspace{57pt} \varepsilon(a_i^{\pm 1}) =1,\hspace{1.2cm} \\
\Ss(u)=u^{n-1}, \qquad  \Ss(x) =-u^{n-1}x, \qquad \Ss(a_i^{\pm 1})=a_i^{\mp 1}, \qquad i,j \in I. \hspace{-1.2cm}
\end{aligned}
\end{equation*}
\end{remark}
\vspace{5pt}

\subsection{Finiteness over the Hopf socle}
The chosen basis $B$ of $\nombre$ \eqref{basis} yields a coalgebra decomposition $\nombre= V_{(0,0)} \oplus \big(\oplus_{(F, E) \in \Gamma} V_{(F,E)}\big)$, where $V_{(F,E)}$ is the subspace spanned by $x^su^ta_F^E$ with $0 \leq s < n,\, 0 \leq t < m$. Observe that $V_{(F,E)}$ is a subcoalgebra by \eqref{comultform}. This decomposition
will be needed to characterize when $\nombre$ is of finite type over its Hopf socle $\dsoc$, being the next result the essential point of the proof in the case $\alpha \neq 0$.

\begin{proposition}\label{block}
Assume that $\alpha \neq 0$.
\begin{enumerate}
\item[(i)] If $q_F^{nE}=1$, then $V_{(F,E)} \simeq T_n(\omega)$ as coalgebras.
\item[(ii)] If $q_F^{nE} \neq 1$, then $V_{(F,E)} \simeq M_n^c(\ku)$ as coalgebras.
\end{enumerate}
\end{proposition}

\pf (i) In view of \eqref{comultform}, $a_F^E$ is group-like in this case. Since the multiplication by a group-like is a coalgebra automorphism, $V_{(F,E)}=T_n(\omega)a_F^E \simeq T_n(\omega)$ as coalgebras. \par \medskip

(ii) For convenience, we abbreviate through the proof $\mu=\alpha(1-q_F^{nE})$ and $a=a_F^E$.
Let $\{c_{st}\}_{1\leq s,t\leq n}$ be the canonical basis of the matrix coalgebra $M_n^c(\ku)$.
We will prove that the map $\Phi:M_n^c(\ku) \rightarrow V_{(F,E)}$ defined by
\begin{equation*}\label{isomorphism}
\xymatrix{c_{st}   \ar@{|->}[0,1]^-{\Phi}
&  {\begin{cases}
{\displaystyle {\binom{s-1}{t-1}}_{\!\omega}\ x^{s-t}u^{t-1}a} & \textrm{if } s \geq t,  \\ \\
{\displaystyle \mu\frac{ (s-1)!_{\omega}}{(t-1)!_{\omega}(n+s-t)!_{\omega}}\ x^{n+s-t}u^{t-1}}a & \textrm{if } s<t,
\end{cases}} }
\end{equation*}
is a coalgebra isomorphism. Clearly, $\Phi$ is bijective and $\varepsilon \Phi(c_{st})=\varepsilon(c_{st})$. To show that $\Phi$ is comultiplicative we must distinguish three cases: \par \vspace{0.4cm}
\noindent (a) Assume that $s=t$:
\renewcommand{\theequation}{\arabic{equation}}
\setcounter{equation}{0}
\begin{align}
\notag
(\Phi \otimes \Phi)\Delta(c_{ss}) & = \sum_{k=1}^{s-1} \Phi(c_{sk}) \otimes \Phi(c_{ks})+\Phi(c_{ss}) \otimes \Phi(c_{ss})+\sum_{k=s+1}^n \Phi(c_{sk}) \otimes \Phi(c_{ks}) \\ \label{eq:2.15}
& \hspace{-0.5cm} = \sum_{k=1}^{s-1} \binom{s-1}{k-1}_{\!\omega} x^{s-k}u^{k-1}a \otimes \mu \frac{(k-1)!_{\omega}}{(s-1)!_{\omega}(n+k-s)!_{\omega}} x^{n+k-s}u^{s-1}a  \\
& \hspace{-0.1cm} + u^{s-1}a \otimes u^{s-1}a \notag \\ \label{eq:2.16}
& \hspace{-0.1cm} + \sum_{k=s+1}^n \mu \frac{(s-1)!_{\omega}}{(k-1)!_{\omega}(n+s-k)!_{\omega}} x^{n+s-k}u^{k-1}a \otimes
\binom{k-1}{s-1}_{\!\omega} x^{k-s}u^{s-1}a
\end{align}
\noindent Put $s-k=v$ in \eqref{eq:2.15} and $k-s=v$ in \eqref{eq:2.16}. Taking into account these changes, we rewrite the $\omega$-coefficients occurring here:
\begin{align*}
\binom{s-1}{k-1}_{\!\omega}\frac{(k-1)!_{\omega}}{(s-1)!_{\omega}(n+k-s)!_{\omega}}
& = \frac{1}{(v)!_{\omega}(n-v)!_{\omega}}.  \\*[5pt]
\frac{(s-1)!_{\omega}}{(k-1)!_{\omega}(n+s-k)!_{\omega}}\binom{k-1}{s-1}_{\!\omega}
&= \frac{1}{(v)!_{\omega}(n-v)!_{\omega}}.
\end{align*}
Substituting all this in the previous equality we have:
\begin{align*}
(\Phi \otimes \Phi)\Delta(c_{ss}) &= u^{s-1}a \otimes u^{s-1}a
+ \mu \sum_{v=1}^{s-1} \frac{1}{(v)!_{\omega}(n-v)!_{\omega}} x^{v}u^{s-v-1}a \otimes x^{n-v}u^{s-1}a \\
& \quad + \mu \sum_{v=1}^{n-s} \frac{1}{(v)!_{\omega}(n-v)!_{\omega}} x^{n-v}u^{v+s-1}a \otimes x^{v}u^{s-1}a \\
& \overset{\text{\ding{172}}}= u^{s-1}a \otimes u^{s-1}a +\mu\es\sum_{k=n-s+1}^{n-1}\es \frac{1}{(k)!_{\omega}(n-k)!_{\omega}} x^{n-k}u^{k+s-1}a\es \otimes\es x^{k}u^{s-1}a  \\
& \quad + \mu \sum_{k=1}^{n-s} \frac{1}{(k)!_{\omega}(n-k)!_{\omega}} x^{n-k}u^{k+s-1}a \otimes x^{k}u^{s-1}a  \\
&= u^{s-1}a \otimes u^{s-1}a + \mu\sum_{k=1}^{n-1} \frac{1}{(k)!_{\omega}(n-k)!_{\omega}} x^{n-k}u^{k+s-1}a \otimes x^{k}u^{s-1}a  \\
&=\Delta(u^{s-1}a) \qquad \textrm{by\ \eqref{comultform}} \\
&= \Delta \Phi(c_{ss}).
\end{align*}

$\text{\ding{172}}$: Put $v=n-k$ in the first sum and $v=k$ in the second.

\medbreak
\setcounter{equation}{0}
\noindent(b) Assume that $s<t$:
\begin{align}
\notag
(\Phi \otimes \Phi)\Delta(c_{st}) &= \sum_{k=1}^{s} \Phi(c_{sk}) \otimes \Phi(c_{kt})+\sum_{k=s+1}^{t-1} \Phi(c_{sk}) \otimes \Phi(c_{kt}) + \sum_{k=t}^n \Phi(c_{sk}) \otimes \Phi(c_{kt})  \notag \\ \label{eq:2.19}
& \hspace{-1cm} = \sum_{k=1}^{s} \binom{s-1}{k-1}_{\!\omega} x^{s-k}u^{k-1}a \otimes \mu \frac{(k-1)!_{\omega}}{(t-1)!_{\omega}(n+k-t)!_{\omega}}
x^{n+k-t}u^{t-1}a \hspace{-20pt} \\ \label{eq:2.20}
& \hspace{-0.5cm} +\sum_{k=s+1}^{t-1} \mu \frac{(s-1)!_{\omega}}{(k-1)!_{\omega}(n+s-k)!_{\omega}} x^{n+s-k}u^{k-1}a \\ \notag
& \otimes \mu \frac{(k-1)!_{\omega}}{(t-1)!_{\omega}(n+k-t)!_{\omega}} x^{n+k-t}u^{t-1}a \\ \label{eq:2.21}
& \hspace{-0.5cm} + \sum_{k=t}^n \mu \frac{(s-1)!_{\omega}}{(k-1)!_{\omega}(n+s-k)!_{\omega}}x^{n+s-k}u^{k-1}a \otimes
\binom{k-1}{t-1}_{\!\omega} x^{k-t}u^{t-1}a
\end{align}
Set $s-k=v$ in \eqref{eq:2.19}, $k-s=v$ in \eqref{eq:2.20}, and $k-t=v$ in \eqref{eq:2.21}. Next we rewrite
the $\omega$-coefficients appearing in the previous equality:

\medbreak
\noindent Coefficient in \eqref{eq:2.19}:
\begin{align*}
 \binom{s-1}{k-1}_{\!\omega} \frac{(k-1)!_{\omega}}{(t-1)!_{\omega}(n+k-t)!_{\omega}}
 &  =  \frac{(s-1)!_{\omega}}{(t-1)!_{\omega}(n+s-t)!_{\omega}} \binom{n+s-t}{v}_{\!\omega}.
\end{align*}
\bigbreak

\noindent Coefficient in \eqref{eq:2.20}:
\begin{align*}
\frac{(s-1)!_{\omega}}{(k-1)!_{\omega}(n+s-k)!_{\omega}} \frac{(k-1)!_{\omega}}{(t-1)!_{\omega}(n+k-t)!_{\omega}} &
\es =\es \frac{(s-1)!_{\omega}}{(t-1)!_{\omega}} \frac{1}{(n-v)!_{\omega}(n+v+s-t)!_{\omega}}.
\end{align*}
\bigbreak

\noindent Coefficient in \eqref{eq:2.21}:
\begin{align*}
\frac{(s-1)!_{\omega}}{(k-1)!_{\omega}(n+s-k)!_{\omega}} \binom{k-1}{t-1}_{\!\omega}
&  = \frac{(s-1)!_{\omega}}{(t-1)!_{\omega}(n+s-t)!_{\omega}} \binom{n+s-t}{v}_{\!\omega}.
\end{align*}
\bigbreak

\noindent Substituting all this in our previous computation we get: \vspace{-8pt}

\begin{align}
& \hspace{0.1cm} = \mu \frac{(s-1)!_{\omega}}{(t-1)!_{\omega}(n+s-t)!_{\omega}} \sum_{v=0}^{s-1} \binom{n+s-t}{v}_{\!\omega}
 x^{v}u^{n-v+s-1}a \otimes x^{n-v+s-t}u^{t-1}a \nonumber \\ \label{eq:2.23}
& \hspace{0.5cm} +\es \mu^2 \frac{(s-1)!_{\omega}}{(t-1)!_{\omega}}\es \sum_{v=1}^{t-s-1}\es \frac{1}{(n-v)!_{\omega}(n+v+s-t)!_{\omega}}\es
x^{n-v}u^{v+s-1}a\es\es \otimes\es x^{n+v+s-t}u^{t-1}a  \\ \label{eq:2.24}
& \hspace{1cm} + \mu \frac{(s-1)!_{\omega}}{(t-1)!_{\omega}(n+s-t)!_{\omega}} \sum_{v=0}^{n-t} \binom{n+s-t}{v}_{\!\omega}
x^{n-v+s-t}u^{v+t-1}a \otimes x^{v}u^{t-1}a
\end{align}
Set $k=n+v+s-t$ in \eqref{eq:2.23} and $l=n-v+s-t$ in \eqref{eq:2.24}. Replacing this in the preceding equality we obtain:
\begin{align}
\label{eq:2.25}
\hspace{-0.3cm} = \mu \frac{(s-1)!_{\omega}}{(t-1)!_{\omega}(n+s-t)!_{\omega}} \sum_{l=0}^{s-1} \binom{n+s-t}{l}_{\!\omega}
x^{l}u^{n-l+s-1}a \otimes x^{n-l+s-t}u^{t-1}a \\
+ \mu^2 \es \frac{(s\es -\es 1)!_{\omega}}{(t\es -\es 1)!_{\omega}}\es \sum_{k=n+s-t+1}^{n-1}\es\es \frac{1}{(k)!_{\omega}(n\es -\es\es k\es\es +\es\es n \es\es + \es s \es -t)!_{\omega}} x^{n-k+n+s-t}u^{k+t-1}a\es\es \otimes \es x^{k}u^{t-1}a \hspace{-0.8cm}  \notag \\ \label{eq:2.28}
+ \mu \frac{(s-1)!_{\omega}}{(t-1)!_{\omega}(n+s-t)!_{\omega}}\es\es \sum_{l=s}^{n+s-t} \es\es \binom{n+s-t}{l}_{\!\omega}\es
x^{l}u^{n-l+s-1}a\es\es \otimes\es x^{n-l+s-t}u^{t-1}a \hspace{-0.5cm}
\end{align}
We join \eqref{eq:2.25} and \eqref{eq:2.28} in a single formula. We now have:
\begin{align*}
& \hspace{0.4cm} = \mu \frac{(s-1)!_{\omega}}{(t-1)!_{\omega}(n+s-t)!_{\omega}}\bigg[ \sum_{l=0}^{n+s-t} \binom{n+s-t}{l}_{\!\omega}
 x^{l}u^{n-l+s-1}a \otimes x^{n-l+s-t}u^{t-1}a  \\
& \hspace{0.7cm} + \mu\es (n\es +\es s\es -\es t)!_{\omega}\es\es\es\es \sum_{k=n+s-t+1}^{n-1} \es\es \frac{1}{(k)!_{\omega}(n\es -\es k\es +\es n\es +\es s\es -\es t)!_{\omega}} x^{n-k+n+s-t}u^{k+t-1}\es a\es\es \otimes\es x^{k}u^{t-1}\es a\bigg]  \\
& \hspace{0.4cm} = \Delta\bigg(\mu \frac{(s-1)!_{\omega}}{(t-1)!_{\omega}(n+s-t)!_{\omega}}x^{n+s-t}u^{t-1}a\bigg)
 \qquad \textrm{by\ \eqref{comultform}}   \\
& \hspace{0.4cm} = \Delta \Phi(c_{st}).
\end{align*}
\setcounter{equation}{0}

\noindent (c) Assume that $s>t$:
\begin{align}
(\Phi \otimes \Phi)\Delta(c_{st}) = \sum_{k=1}^{t-1} \Phi(c_{sk}) \otimes \Phi(c_{kt}) + \sum_{k=t}^{s} \Phi(c_{sk}) \otimes \Phi(c_{kt}) + \sum_{k=s+1}^n \Phi(c_{sk}) \otimes \Phi(c_{kt}) \notag \hspace{1.9cm} \\[-2pt]
= \sum_{k=1}^{t-1} \binom{s-1}{k-1}_{\!\omega} x^{s-k}u^{k-1}a \otimes \mu \frac{(k-1)!_{\omega}}{(t-1)!_{\omega}(n+k-t)!_{\omega}} x^{n+k-t}u^{t-1}a \label{eq:2.29} \hspace{2.3cm} \\ \label{eq:2.30}
+ \sum_{k=t}^{s} \binom{s-1}{k-1}_{\!\omega} x^{s-k}u^{k-1}a \otimes \binom{k-1}{t-1}_{\!\omega}
x^{k-t}u^{t-1}a \hspace{4.5cm} \\ \label{eq:2.31}
+ \sum_{k=s+1}^n \mu \frac{(s-1)!_{\omega}}{(k-1)!_{\omega}(n+s-k)!_{\omega}} x^{n+s-k}u^{k-1}a \otimes
\binom{k-1}{t-1}_{\!\omega} x^{k-t}u^{t-1}a. \hspace{1.7cm}
\end{align}

Put $v=n+k-t$ in \eqref{eq:2.29}, $l=s-k$ in \eqref{eq:2.30}, and $v=k-t$ in \eqref{eq:2.31}. We rewrite the $\omega$-coefficients occurring here
taking into account these changes:
\medbreak

\noindent Coefficient in \eqref{eq:2.29}:
\begin{align*}
\binom{s-1}{k-1}_{\!\omega}\frac{(k-1)!_{\omega}}{(t-1)!_{\omega}(n+k-t)!_{\omega}}
&  = \binom{s-1}{t-1}_{\!\omega}\ \frac{(s-t)!_{\omega}}{(v)!_{\omega}(n-v+s-t)!_{\omega}}.
\end{align*}
\medbreak

\noindent Coefficient in \eqref{eq:2.30}:
\begin{align*}
\binom{s-1}{k-1}_{\!\omega} \binom{k-1}{t-1}_{\!\omega} &
= \binom{s-1}{t-1}_{\!\omega} \binom{s-t}{l}_{\!\omega}.
\end{align*}
\medbreak

\noindent Coefficient in \eqref{eq:2.31}:
\begin{align*}
\frac{(s-1)!_{\omega}}{(k-1)!_{\omega}(n+s-k)!_{\omega}} \binom{k-1}{t-1}_{\!\omega} &
= \binom{s-1}{t-1}_{\!\omega}\ \frac{(s-t)!_{\omega}}{(v)!_{\omega}(n-v+s-t)!_{\omega}}.
\end{align*}
\bigbreak

\noindent Substituting all this in our previous computation we have:
\begin{align}
= \mu \binom{s-1}{t-1}_{\!\omega}\es\es (s-t)!_{\omega}\es\es \sum_{v=n-t+1}^{n-1}\es \frac{1}{(v)!_{\omega}(n-v+s-t)!_{\omega}}
x^{n-v+s-t}u^{v+t-1}a\es\es \otimes\es x^{v}u^{t-1}a \label{eq:2.32} \\
+ \binom{s-1}{t-1}_{\!\omega} \sum_{l=0}^{s-t} \binom{s-t}{l}_{\!\omega} x^lu^{s-t-l+t-1}a \otimes x^{s-t-l}u^{t-1}a \hspace{3.75cm} \notag\\
+ \mu \binom{s\es -\es 1}{t\es -\es 1}_{\!\omega}\es\es (s\es -\es t)!_{\omega}\es\es\es \sum_{v=s-t+1}^{n-t}\es\es \frac{1}{(v)!_{\omega}(n\es -\es v\es +\es s\es -\es t)!_{\omega}}
x^{n-v+s-t}u^{v+t-1}a\es\es \otimes\es x^{v}u^{t-1}a \label{eq:2.36}
\end{align}
We join \eqref{eq:2.32} and \eqref{eq:2.36} in a single formula. Then we obtain:
\begin{align*}
&= \binom{s-1}{t-1}_{\!\omega}\left[\sum_{l=0}^{s-t} \binom{s-t}{l}_{\!\omega} x^lu^{s-t-l+t-1}a \otimes x^{s-t-l}u^{t-1}
a \right.   \\
& \qquad + \mu \left. (s-t)!_{\omega} \sum_{v=s-t+1}^{n-1} \frac{1}{(v)!_{\omega}(n-v+s-t)!_{\omega}} x^{n-v+s-t}u^{v+t-1}a \otimes x^{v}u^{t-1}a\right] \\
& = \Delta\left(\binom{s-1}{t-1}_{\!\omega} x^{s-t}u^{t-1}a\right) \qquad \textrm{by\ \eqref{comultform}} \\
& = \Delta\Phi(c_{st}).
\end{align*}
\epf

We are now ready to characterize when $\nombre$ is of finite type over $\dsoc$, from which Theorem \ref{question-ad} in the introduction will follow.

\setcounter{equation}{14}

\begin{theorem}\label{mainth2}
The Hopf algebra $\nombre$ is co-Frobenius and $\dsoc=\ku G(\nombre)$. Moreover:
\begin{enumerate}
\item[(i)] If $\alpha=0$, then $\nombre$ is of finite type over $\dsoc.$ \vspace{3pt}
\item[(ii)] If $\alpha \neq 0$, then $\nombre$ is of finite type over $\dsoc$ if and only if there is a finite subset $J$ of $I$ such that $q_i^n=1$ for all $i \in I\backslash J$ and $q_j$ is an $\nu_j$-th root of unity for all $j \in J$.
\end{enumerate}
\end{theorem}

\pf Let $\mathcal{A}$ be the (finite dimensional) Hopf subalgebra of $\nombre$ generated by $u$ and $x$. The particular form of the chosen basis $B$  of $\nombre$ \eqref{basis}, together with \eqref{comultform}, gives us a coalgebra decomposition $\nombre= \mathcal{A} \oplus \big(\oplus_{(F,E) \in \Gamma} V_{\!(F,E)}\big)$. Viewed as a right $\nombre$-comodule, $\mathcal{A}$ is injective and contains $\ku$. Hence $E(\ku) \subset \mathcal{A}$, so that $\nombre$ is co-Frobenius. \par \smallskip

We next describe $\dsoc$ and show that it equals $\ku G(\nombre)$. The coalgebra decomposition above is inherited to the coradical, that is, $\nombre_0=\ku C_m \oplus \big(\oplus_{(F,E) \in \Gamma} {V_{(F,E)}}_0 \big)$, where $C_m$ is generated by $u$. Then,
\begin{align}\label{hopfsocledec}
\dsoc & = {\displaystyle (\ku C_m \cap \dsoc) \oplus \big(\oplus_{(F,E) \in \Gamma} ({V_{(F,E)}}_0 \cap \dsoc )\big)} \nonumber \\
 &  = {\displaystyle \ku C_m \oplus \big(\oplus_{(F,E) \in \Gamma} ({V_{(F,E)}}_0 \cap \dsoc) \big).}
\end{align}
Here we used that $G(\nombre)\subseteq \dsoc.$ \par \bigskip

(i) Assume that $\alpha=0$. Then $u^ta_F^E \in G(\nombre)$ for all $(F,E) \in \Gamma \cup \{(0,0)\}$ and $0 \leq t < m.$ Setting $X=\{x^s: 0 \leq s <n\}$, we have $\nombre =X \ku G(\nombre) \subseteq X\dsoc.$ Thus $\nombre$ is of finite type over $\dsoc$. To see that $\dsoc =\ku G(\nombre)$, notice that the multiplication by a group-like element establishes a coalgebra automorphism of $\nombre$. Then $V_{\!(F,E)}=\mathcal{A}a_F^E \simeq \mathcal{A}$ as coalgebras. This implies ${V_{\!(F,E)}}_0=\mathcal{A}_0a_F^E=\oplus_{t=0}^{m-1} \ku u^ta_F^E$. Consequently, ${V_{\!(F,E)}}_0 \cap \dsoc \subset \ku G(\nombre)$. \par \bigskip

(ii) Assume that $\alpha \neq 0$. Recall that in this case $n=m$ and $\mathcal{A}=T_n(\omega)$. To compute $\dsoc$ we first calculate ${V_{\!(F,E)}}_0$ and see if it is contained or not there. We must distinguish two cases: $q_F^{nE} \neq 1$ and $q_F^{nE}=1.$ \par \smallskip

{\it Case} $q_F^{nE} \neq 1.$ By Proposition \ref{block} (ii), $V_{\!(F,E)} \simeq M_n^c(\ku)$ as coalgebras. Then $V_{\!(F,E)}={V_{\!(F,E)}}_0$. Let $S_{(F,E)}$ be the unique (up to isomorphism) simple right $\nombre$-comodule corresponding to $V_{\!(F,E)}.$ From the aforementioned coalgebra decomposition of $\nombre$, we conclude that $S_{(F,E)}$ is injective. We claim that $S_{(F,E)} \nsubset \dsoc$. Otherwise, $S_{(F,E)} \otimes S_{(F,E)}^*$ would be semisimple and injective. As it contains $\ku$, it would follow that $\ku$ is injective and therefore $\nombre$ would be cosemisimple. This is not possible because $T_n(\omega)$ is a non cosemisimple Hopf subalgebra of $\nombre$. \par \smallskip

{\it Case} $q_F^{nE}=1.$ Here we argue as for $\alpha=0$. By \eqref{comultform}, $a_F^E \in G(\nombre).$ By Proposition \ref{block} (i), $V_{\!(F,E)} \simeq T_n(\omega)$ as coalgebras and, consequently, ${V_{\!(F,E)}}_0$ is spanned by $u^ta_F^E$ with $t=0,\dots,n-1$. This implies ${V_{\!(F,E)}}_0 \subset \dsoc.$
\par \smallskip

Let $\Lambda =\{(F,E) \in \Gamma: q_F^{nE}=1\} \cup \{(0,0)\}$ and $\bar{\Lambda}=\{(F,E) \in \Gamma: q_F^{nE}\neq 1\}.$ Our previous discussion, together with the decomposition (\ref{hopfsocledec}) of $\dsoc$, entails $\dsoc = \oplus_{t=0}^{n-1} \oplus_{(F,E) \in \Lambda} \ku u^ta_F^E \subset \ku G(\nombre)$. Every $g \in G(\nombre)$ must be then of the form $u^ta_F^E$ for some $t$ and $(F,E) \in \Lambda$. We write $R_{t,(F,E)}$ for the simple right $\nombre$-comodule
$\ku u^ta_F^E$. \par \medskip

We next proceed to prove the statement about the finite generation over $\dsoc$. \par \smallskip

Suppose that there is $J \subseteq I$ finite such that $q_i^n=1$ for all $i \in I\backslash J$ and $q_j$ is an $\nu_j$-th root of unity for all
$j \in J$. Since $q_j^{\nu_j}=1$, in view of \eqref{comultform}, $a_{j}^{\nu_j} \in G(\nombre)$ for all $j \in J$. Similarly, $a_i \in G(\nombre)$ for all $i \in I\backslash J$. Put $J=\{j_1,\dots,j_l\}$. Consider the set $$X=\{x^sa_{j_1}^{f_{j_1}}\dots\hspace{1pt} a_{j_l}^{f_{j_l}}: 0 \leq s< n, 0 \leq f_{j_k} <\nu_{j_k}, 1 \leq k \leq l\}.$$ We prove that $\nombre=X\dsoc$. It suffices to show the inclusion for elements of the form $a_{F}^E$. Set $F=(i_1,\dots,i_r)$ and $E=(e_1,\dots,e_r)$, so that $a_{F}^E=a_{i_1}^{e_1}\dots\hspace{1pt} a_{i_r}^{e_r}$. We can assume that $i_k \in J$ for $k=1,\dots,p$ and $i_k \notin J$ for $k=p+1,\dots,r$. Then $a_{i_k} \in G(\nombre)$ for $k=p+1,\dots,r.$ Write $e_k=\nu_kc_k+\bar{e}_k$ with $0 \leq \bar{e}_k < \nu_k$ for $k=1,\dots,p$. We now have
$a_{F}^E= (a_{i_1}^{\bar{e}_1}\dots a_{i_p}^{\bar{e}_p})(a_{i_1}^{\nu_1 c_1} \dots a_{i_p}^{\nu_p c_p}a_{i_{p+1}}^{e_{p+1}} \dots a_{i_r}^{e_r}) \in XG(\nombre).$
\par \smallskip

Conversely, assume that $\nombre$ is of finite type over $\dsoc$. Take a finite set $Y \subset \nombre$ such that $\nombre=Y\dsoc$; then $\nombre=W \dsoc$ for the (finite dimensional) subcoalgebra $W$ generated by $Y$. We have a coalgebra decomposition
$$W= (W \cap \mathcal{A}) \oplus \big(\oplus_{(F,E) \in \Gamma} (W \cap V_{\! (F,E)}) \big)$$ inherited from that of $\nombre$, see \cite[Exercise 2.2.18 (iv)]{DNR} and the proof of \cite[Lemma 5.1.9]{M}. As $\dim W < \infty$, the set $\Omega=\{(F,E) \in \Gamma : W \cap V_{\!(F,E)}\neq 0\}$ is finite. We can express $\Omega$ as a disjoint union of $\Omega_1=\{(F,E) \in \Omega: q_F^{nE}=1\}$ and $\Omega_2=\{(F,E) \in \Omega: q_F^{nE} \neq 1\}$. Then $a_F^E \in G(\nombre)$ for all $(F,E) \in \Omega_1$. For $(F,E) \in \Omega_2$ set $F=(i_1,\dots,i_r)$ and $E=(e_1,\dots,e_r)$, so that $a_{F}^E=a_{i_1}^{e_1} \dots a_{i_r}^{e_r}$. Remove from $F$ those $i_k$ such that $q_{i_k}^n=1\ (a_{i_k} \in G(\nombre)).$ In this way,  $\Omega_2$ gives rise to a set $J=\{j_1,\dots,j_l\}$ such that $q_{j_k}^n \neq 1$ for all $k=1,\dots,l$. \par \smallskip

Since $\dsoc=\ku G(\nombre)$, we have another coalgebra decomposition $\nombre=\oplus_{g \in G(\nombre)} Wg.$ Let $i \in I$ be such that $q_i^n \neq 1$ and consider the simple subcoalgebra $\mathcal{A}a_i$, Proposition \ref{block} (ii). There is $g_{i} \in G(\nombre)$ such that $\mathcal{A}a_i \subset Wg_i$. We know that $g_i^{-1}=u^ta_{F}^E$ for certain $(F,E)\in \Omega_1$. Then  $a_ig_{i}^{-1}$ is of the form $u^ta_{F'}^{E'}$ with $q_{F'}^{nE'} \neq 1$ and $V_{\! (F',E')}=\mathcal{A}a_ig_{i}^{-1} \subset W$. Hence $(F',E') \in \Omega_2$ and so $i \in J$. This shows that $i \in J$ if and only if $q_i^n \neq 1$. \par \smallskip

We must finally prove that $q_{j_k}$ is an $\nu_k$-th root of unity for every $k \in \{1,\dots,l\}$.
This is clear if the ground field $\ku$ is finite. Assume that $\ku$ is infinite. Let $e_{j_k}$ be the maximum exponent, in absolute value, of $a_{j_k}$ when occurring in the elements $a_F^E$ with $(F,E) \in \Omega_2$. Put
$$X=\{x^sa_{j_1}^{\pm f_{j_1}}\dots\ a_{j_l}^{\pm f_{j_l}}: 0 \leq s< n, 0 \leq f_{j_k} \leq e_{j_k}, 1 \leq k \leq l\}.$$
Then $\nombre=W\dsoc \subseteq X\dsoc.$ Set $P=\{\pm f_{j_k} : 1 \leq k \leq l\}.$ Pick $z \in \Ent$ and $k \in \{1,\dots,l\}$ arbitraries. We write
$$a_{j_k}^z = \sum_{\sigma} \lambda_{\sigma} y_{\sigma} a_{F_{\sigma}}^{E_{\sigma}}, \qquad \  \lambda_{\sigma} \in \ku, \ y_{\sigma}  \in X, \ a_{F_{\sigma}}^{E_{\sigma}} \in G(\nombre).$$
Observe that $q_{F_{\sigma}}^{nE_{\sigma}}=1$ for all $\sigma$.
We assume that all terms $y_{\sigma} a_{F_{\sigma}}^{E_{\sigma}}$ in this sum \vspace{-2pt} are distinct.
A priori, $a_{F_{\sigma}}^{E_{\sigma}}$ might be multiplied by a power of $u$, as $u \in G(\nombre)$, but by linear independence,
this is not possible. By the same reason, no power of $x$ occurs in $y_{\sigma}.$
Set $F_{\sigma}=(i_{1},\dots,i_{r_{\sigma}}), E_{\sigma}=(\theta_{\sigma,1},\dots,\theta_{\sigma, {r_{\sigma}}})$ and
\begin{equation*}
y_{\sigma} = a_{j_1}^{\gamma_{\sigma,1}}\dots\ a_{j_l}^{\gamma_{\sigma,l}}, \qquad \gamma_{\sigma,1},\dots, \gamma_{\sigma,l} \in P.
\end{equation*}
Set $F'_{\sigma}=F_{\sigma}\backslash (F_{\sigma} \cap J)$. Take out of $a_{F_{\sigma}}^{E_{\sigma}}$ the $a_i$'s with $i \in F_{\sigma} \cap J$, join them to $y_{\sigma}$ and consider the corresponding list of exponents $E'_{\sigma}$ for $F'_{\sigma}$. Then,
\begin{equation}\label{sumel2}
a_{j_k}^z= \sum_{\sigma} \lambda_{\sigma}a_{j_1}^{\gamma_{\sigma,1}+\theta_{\sigma,1}}\dots\ a_{j_l}^{\gamma_{\sigma,l}+\theta_{\sigma,l}}a_{F'_{\sigma}}^{E'_{\sigma}}.
\end{equation}
Here we are abusing a bit of notation because some $j_{\mu}$ could not appear in $F_{\sigma}$ and, in such a case, we understand $\theta_{\sigma,\mu}=0$.
Moreover, $\theta_{\sigma,1},\dots,\theta_{\sigma,l}$ correspond to $a_{i_1},\dots,a_{i_l}$ and subindices $i_1,\dots,i_l$ are not necessarily equal $j_1,\dots,j_l$.
All monomials in the right-hand side of \eqref{sumel2} are different basis element. There must \vspace{-1pt} be $\tau$ such that $a_{j_k}^z=a_{j_k}^{\gamma_{\tau,k}+\theta_{\tau,k}}$ and $a_{j_{\mu}}^{\gamma_{\tau,\mu}+\theta_{\tau,\mu}}=1$ for $\mu \neq k$. From here, $z=\gamma_{\tau,k}+\theta_{\tau,k}$ and\vspace{-2pt} $\theta_{\tau,\mu}=-\gamma_{\tau,\mu}$. Recall that $q_{F_{\tau}}^{nE_{\tau}}=q_{i_1}^{n\theta_{\tau,1}}\dots\ q_{i_{r_{\tau}}}^{n\theta_{\tau, r_{\tau}}}=1$ and $q_i^n=1$ for $i \notin J$.\vspace{-2pt} Then  $q_{j_1}^{\theta_{\tau,1}}\dots\ q_{j_{l}}^{\theta_{\tau,l}}$ is an $n$-th root of unity, say $\zeta$. Substituting the value of the $\theta$'s just found we get:
$$q_{j_k}^z= \zeta \prod_{\mu=1}^l q_{j_{\mu}}^{\gamma_{\tau,\mu}}.$$
The right-hand side of this equality only takes a finite number of values because the $\gamma$'s are chosen from the finite set $P$ and $\zeta$ is an $n$-th root of unity. However, $z$ runs over $\Ent$. This yields that $q_{j_k}$ is a root of unity.
\epf \bigbreak

As a consequence of Theorem \ref{mainth2} (ii), if some $q_i$ is not a root of unity, then $\nombre$ is not of finite type over $\dsoc$, establishing so Theorem \ref{question-ad} announced in the introduction. On the other hand, if $\alpha \neq 0, \vert I \vert=1$ and $q_i$ is not a root of unity, then $\dsoc$ is finite dimensional. Its only elements are the group-like elements of Taft Hopf algebra. \medskip

The proof of Theorem \ref{mainth2} (ii) yields that when $\alpha \neq 0$ the set
$$\{R_{t,(F,E)}: 0 \leq t < n, (F,E) \in \Lambda\} \cup \{S_{(F,E)}:(F,E) \in \bar{\Lambda}\}$$
is a full set of representatives of the simple right $\nombre$-comodules. Moreover, the ${S_{(F,E)}}$'s are injective.

\begin{proposition}
The multiplication rules for the above set are:
\begin{eqnarray*}
R_{t,(F,E)} \otimes R_{t',(F',E')} \simeq R_{t+t',(F\cup F',E+E')}, \hspace{3.9cm} \vspace{8pt} \\
R_{t,(F,E)} \otimes S_{(F',E')} \simeq S_{(F\cup F',E+E')} \simeq S_{(F',E')} \otimes R_{t,(F,E)}, \vspace{8pt} \hspace{1.5cm} \\ \quad S_{(F,E)} \otimes S_{(F',E')} \simeq \left\{\begin{array}{ll}
T_n(\omega)  & \ {\rm if} \ F' = F\ {\rm and} \ E'=-E, \vspace{2pt} \\
S_{(F\cup F',E+E')}^{\,n} & \ {\rm otherwise}.
\end{array}\right.
\end{eqnarray*}

\pf The first isomorphism is clear because $u^ta_F^E$ and $u^{t'}\es\es a_{F'}^{E'}$ are group-like. The coefficient space of $R_{t,(F,E)} \otimes S_{(F',E')}$ equals $$(\ku u^ta_F^E)V_{(F',E')}=(\ku u^ta_F^E)(T_n(\omega)a_{F'}^{E'})=V_{(F\cup F',E+E')}.$$ Since the latter is a simple subcoalgebra, $R_{t,(F,E)} \otimes S_{(F',E')}$ is isomorphic to $r$ copies of $S_{(F\cup F',E+E')}.$ Comparing dimensions, $r=1.$ Similarly, $S_{(F',E')} \otimes R_{t,(F,E)} \simeq S_{(F\cup F',E+E')}$. Assume finally that either $F' \neq F$ or $E' \neq -E$. Then the coefficient space of $S_{(F,E)} \otimes S_{(F',E')}$ equals $V_{(F,E)}V_{(F',E')}=V_{(F\cup F',E+E')}.$ Now argue as before. In the case $F' = F$ and $E'=-E$, the coefficient space is $T_n(\omega)$. Observe that $S_{(F',E')} \simeq S_{(F,E)}^{\,*}$ because $\Ss(V_{(F,E)})=V_{(F,-E)}$.\vspace{-0.4pt} Then $S_{(F,E)} \otimes S_{(F',E')}$ contains $\ku$. Tensoring with $\ku u^t$ and using the previous isomorphism, $S_{(F,E)} \otimes S_{(F',E')}$ \vspace{-1pt} contains $\ku u^t$ for $t=0,\ldots,n-1$. Since it is injective, it must contain a copy of the injective hull of $\ku u^t$ for $t=0,\ldots,n-1$. These are $n$-dimensional. Comparing dimensions, it must contain exactly one copy of each. Hence $S_{(F,E)} \otimes S_{(F',E')} \simeq T_n(\omega).$ \epf

\end{proposition}

\section{Examples in positive characteristic}\label{sec:positive-car}

In this final section we construct from group theory examples of co-Frobenius Hopf algebras, over fields of positive characteristic, that are not of finite type over their Hopf socles. \par \bigskip

Let $G$ be an infinite group and let $K$ be a finite abelian group of order $n$ acting freely on $G$ by group automorphisms.
Then the group algebra $\ku G$ is a left comodule Hopf algebra over the dual group algebra $\ku^K$ with coaction $\rho:\ku G \rightarrow \ku^K \otimes \ku G$, $g \mapsto \sum_{k \in K} \delta_{k} \otimes (k \cdot g)$. Let $H=\ku G \# \ku^K$ be the smash coproduct Hopf algebra; this is the tensor product algebra, with  comultiplication and antipode given by
\begin{align}\label{comultsmash}
\Delta(g \#  \delta_{k}) & = \sum_{t \in K} (g \#  \delta_{t}) \otimes (t \cdot g \# \delta_{t^{-1}k}),
& \Ss(g \#  \delta_{k}) &= k \cdot g^{-1} \# \delta_{k^{-1}}.
\end{align}

Suppose now that $\car \ku$ divides $n$. Then $\ku^K$ is co-Frobenius but not cosemisimple. Being a cleft extension of co-Frobenius Hopf algebras, $H$ is co-Frobenius by \cite[Proposition 5.2]{BDGN}-- see also \cite[Theorem 2.10]{AC}-- and not cosemisimple because the Hopf subalgebra $\ku^K$ is not so.

\begin{proposition}\label{decomp}
Let $\Gamma$ be a set of representatives of the orbits in $G$. For $g \in G$ let $\bar{g} \in \Gamma$ denote the representative of $\Oc(g)$. Let $\ku \Oc(g) \subset \ku G$ be the $\ku$-vector subspace spanned by $\Oc(g)$. Then:
\begin{enumerate}
\item[(i)] $H=\oplus_{\bar{g} \in \Gamma} \ \ku \Oc(\bar{g}) \# \ku^K$ as coalgebras. \vspace{3pt}
\item[(ii)] $\ku \Oc(\bar{g}) \# \ku^K \simeq M^c_n(\ku)$ as coalgebras for $\bar{g} \neq 1_G$ and $\ku 1_G \# \ku^K\simeq \ku^K$.
\end{enumerate}
\end{proposition}

\pf (i) By \eqref{comultsmash}, $\ku \Oc(\bar{g}) \# \ku^K$ is a subcoalgebra, and clearly $H=\oplus_{\bar{g} \in \Gamma} \ \ku \Oc(\bar{g}) \# \ku^K$. \smallbreak

(ii) Let $\{c_{ij}\}_{1\le i,j\le n}$ be the canonical basis of $M^c_n(\ku)$. We write $K=\{k_1,\dots,k_n\},$ where $k_1=1_K.$ It is not difficult to show that the map $\Phi:\ku \Oc(\bar{g}) \# \ku^K \rightarrow M^c_n(\ku)$, $(k_i \cdot g) \# \delta_{k_j} \mapsto e_{il}$, with $k_ik_j=k_l$, is a coalgebra isomorphism. That $\ku 1_G \# \ku^K \simeq \ku^K$ is clear. \epf

With notation as above, the unique (up to isomorphism) simple right $M^c_n(\ku)$-comodule is $\ku \{c_{1j}: j=1,\dots,n\}$. Then, for $\bar{g} \neq 1_G$, through $\Phi$, the simple right $H$-comodule corresponding to the block $\ku \Oc(g) \# \ku^K$ is $S_{\bar{g}}=\ku\{g \# \delta_{k_j}: j=1,\dots,n\}$. Put $S_{\overline{1_G}}=\ku \{1_G \# \delta_{k_j}: j=1,\dots,n\},$ which is isomorphic to $\ku^K$ as a right $H$-comodule. Notice that this is never simple if $K$ is non trivial.

\begin{theorem}\label{mainth3}
The Hopf algebra $H$ is co-Frobenius and not of finite type over $\hsoc$.
\end{theorem}

\pf We show that $S_{\bar{g}}$ is not included in $\hsoc$ for $\bar{g} \neq 1_G$. From Proposition \ref{decomp} (i), each $S_{\bar{g}}$ is injective. If it is contained in $\hsoc$, then $S_{\bar{g}} \otimes S_{\bar{g}}^*$ is semisimple. Since $S_{\bar{g}} \otimes S_{\bar{g}}^*$
must contain $\ku$ as a direct summand, we would have that $\ku$ is injective and hence $H$ would be cosemisimple, a contradiction. Taking into account again the coalgebra decomposition of $H$ in Proposition \ref{decomp}, only the simple comodules of $\ku^K$ could be included in $\hsoc$. There is a finite number of them (up to isomorphism), so that $\hsoc$ is finite dimensional. Since $G$ is infinite, $H$ cannot be of finite type over $\hsoc$. \epf

We just gave an indirect argument to show that $S_{\bar{g}}$ is not included in $\hsoc$ for $\bar{g} \neq 1_G$. However, it is possible to compute the decomposition of $S_{\bar{g}} \otimes S_{\bar{h}}$:

\begin{proposition}\label{posdecomp}
As right $H$-comodules, $S_{\bar{g}} \otimes S_{\bar{h}} \simeq \oplus_{r=1}^n S_{\overline{(k_r\cdot g)h}}$. In particular, if $\bar{g},\bar{h} \neq 1_G$, then $S_{\bar{g}} \otimes S_{\bar{h}}$ is semisimple except when $\Oc(g)=\Oc(h^{-1})$. Equivalently, when either $S_{\bar{g}} \simeq S_{\bar{h}}^*$ or $S_{\bar{h}} \simeq S_{\bar{g}}^*$.
\end{proposition}

\pf Define $f:S_{\bar{g}} \otimes S_{\bar{h}} \rightarrow \oplus_{r=1}^n S_{\overline{(k_r\cdot g)h}},(g \# \delta_{k_l}) \otimes (h \# \delta_{k_m}) \mapsto (k_r\cdot g)h \# \delta_{k_lk_r^{-1}}$, where $k_r=k_m^{-1}k_l$. Clearly, $f$ is bijective. We check that it is a comodule morphism:
\begin{align*}
(f \otimes \id)&\rho[(g \# \delta_{k_l}) \otimes (h \# \delta_{m})] = \\
&= (f \otimes \id)\bigg[\sum_{i,j=1}^n (g \# \delta_{k_i}) \otimes (h \# \delta_{k_j}) \otimes
 (k_i \cdot g)(k_j \cdot h) \# \delta_{k_i^{-1}k_l}\delta_{k_j^{-1}k_m}\bigg]
\\
&= (f \otimes \id)\bigg[\sum_{i=1}^n (g \# \delta_{k_i}) \otimes (h \# \delta_{k_ik_mk_l^{-1}}) \otimes
 (k_i \cdot g)((k_ik_mk_l^{-1}) \cdot h) \# \delta_{k_i^{-1}k_l}\bigg] \\
& =   \sum_{i=1}^n ((k_m^{-1}k_l)\cdot g)h \# \delta_{k_ik_mk_l^{-1}} \otimes [k_i \cdot (g ((k_mk_l^{-1}) \cdot h))]
\# \delta_{k_i^{-1}k_l} \\
& =   \sum_{j=1}^n ((k_m^{-1}k_l)\cdot g)h \# \delta_{k_j} \otimes [k_j \cdot (((k_m^{-1}k_l) \cdot g)h)] \# \delta_{k_j^{-1}k_m} \\
& =  \rho (((k_m^{-1}k_l) \cdot g)h \# \delta_{k_m})  \\
& =  \rho f[(g \# \delta_{k_l}) \otimes (h \# \delta_{k_m})].
\end{align*}
The comodule $S_{\overline{(k_r\cdot g)h}}$ is simple except when $\overline{(k_r\cdot g)h}=1_G$. This happens if and only if $\Oc(g)=\Oc(h^{-1})$. In this case, $S_{\overline{(k_r\cdot g)h}}=S_{\overline{1_G}}=\ku^K,$ which is not semisimple. Finally, from \eqref{comultsmash} and Proposition \ref{decomp} (ii), it follows that $S_{\overline{h}}^* \simeq S_{\overline{h^{-1}}}$ for $\overline{h} \neq 1_G$. \epf

We finish this paper by providing an example of a commutative infinite dimensional co-Frobenius Hopf algebra whose Hopf socle is trivial.

\begin{example}
Assume that $\car \ku=2$. Let $G=\Ent$ and $K=C_2$ with generators $x$ and $\sigma$ respectively. Let $K$ act on $G$ by $\sigma\cdot x=x^{-1}$. We have $\Oc(1)=\{1\}$ and $\Oc(x^z)=\{x^z,x^{-z}\}$ for $0 \neq z \in \Ent$. We take $\{x^n : n \geq 0\}$ as a set of representatives of the orbits. We construct the smash coproduct Hopf algebra $H=\ku\Ent \# \ku^{C_2}$. The $\ku$-vector subspace $D_n$ spanned by $\{x^{\pm n} \# \delta_{\sigma^i} : i=0,1\}$ is a subcoalgebra of $H$ isomorphic to $M_2^c(\ku)$. By Proposition \ref{decomp}, we have a coalgebra decomposition $H=\ku^{C_2} \oplus (\oplus_{n \in \Na} D_n)$. The simple comodule $S_n$ attached to $D_n$ is spanned by $\{x^{n} \# \delta_{\sigma^i} : i=0,1\}$. The only simple comodule of $\ku^{C_2}$ is $\ku$. Hence $\hsoc=\ku$ by the proof of Theorem \ref{mainth3}. In view of Proposition \ref{posdecomp}, the decomposition rules for the tensor product of the simple comodules are: $S_n \otimes S_m \simeq S_{n+m} \oplus S_{\vert m-n\vert}$ for $m \neq n$ and $S_n \otimes S_n \simeq S_{2n} \oplus \ku^{C_2}$. \smallbreak
\end{example}

\subsection*{Acknowledgements} N. A. was partially supported by CONICET,
FONCyT and Secyt (UNC). J. C. is supported by projects MTM2011-27090 from MICINN and FEDER and P07-FQM03128 from Junta de Andaluc\'{\i}a. The work of P. E. was partially supported by the NSF grant DMS-1000113.

%\newpage


\begin{thebibliography}{WWA}

\bibitem[A]{A}  H. H. Andersen, \emph{Tensor products of quantized tilting
modules}. Comm. Math. Phys. {\bf 149} (1992), 149--159.

\bibitem[APW]{APW} H. H. Andersen, P. Polo and K. Wen, \emph{Representations of quantum algebras}. Invent. Math. {\bf 104} (1991), 1-59.

\bibitem[AC]{AC} N. Andruskiewitsch and J. Cuadra, \emph{On the structure of (co-Frobenius) Hopf algebras}. J. Noncommut. Geom. {\bf 7} (2013), 83-104.

\bibitem[AD]{AD} N. Andruskiewitsch and S. D\u{a}sc\u{a}lescu, \emph{Co-Frobenius Hopf algebras and the coradical filtration}. Math.
Z. \textbf{243} (2003), 145-154.

\bibitem[AS]{AS1} {N. Andruskiewitsch} and  {H.-J. Schneider}, \emph{Lifting of quantum linear spaces and pointed Hopf algebras of order $p{^3}$}. J. Algebra \textbf{209} (1998), 658-691.

\bibitem[BDGN]{BDGN} \hspace{-6pt} M. Beattie, S. D\u{a}sc\u{a}lescu, L. Gr\"{u}nenfelder and C. N\u{a}st\u{a}sescu, \emph{Finiteness
conditions, co-Frobenius Hopf algebras, and quantum groups}. J. Algebra {\bf 200} (1998), 312-333.

\bibitem[C]{C} J. Cuadra, \emph{On Hopf algebras with nonzero integral}. Comm. Algebra \textbf{34} (2006), 2143-2156.

\bibitem[De1]{De1} P. Deligne, \emph{Catégories tensorielles.} Mosc. Math. J. {\bf 2} (2002), 227-248.

\bibitem[De2]{De2} \bysame, \emph{La catégorie des représentations du groupe symétrique $S_t$, lorsque $t$ n'est pas un entier naturel.} Algebraic groups and homogeneous spaces, 209-273. Tata Inst. Fund. Res. Stud. Math., Tata Inst. Fund. Res., Mumbai, 2007.

\bibitem[DN]{DN} S. D\u{a}sc\u{a}lescu and C. N\u{a}st\u{a}sescu, \emph{Coactions on spaces of morphisms.} Algebr. Represent. Theory {\bf 12} (2009), 193-198.

\bibitem[DNR]{DNR} S. D\u{a}sc\u{a}lescu, C. N\u{a}st\u{a}sescu, and \c{S}. Raianu, \emph{Hopf algebras. An introduction}. Monographs and Textbooks in Pure and Applied Mathematics {\bf 235}. Marcel-Dekker, New-York, 2001.

\bibitem[Do1]{Don1} S. Donkin, {\it On projective modules for algebraic groups.} J. London Math. Soc. (2) {\bf 54} (1996), 75-88.

\bibitem[Do2]{Don2} S. Donkin, \emph{On the existence of Auslander-Reiten sequences of group representations II}. Algebr. Represent. Theory
{\bf 1} (1998), 215-253.

\bibitem[EGNO]{EGNO} \hspace{-6pt} P. Etingof, S. Gelaki, D. Nikshych, and V. Ostrik, \emph{Tensor categories}. Lecture notes available at http://euclid.unh.edu/\~{}nikshych/.

\bibitem[EO]{EO} P. Etingof, V. Ostrik, \emph{Finite tensor categories}. Mosc. Math. J. {\bf 4} (2004), 627-654.

\bibitem[G]{G} J. A. Green, \emph{Locally finite representations}. J. Algebra {\bf 41} (1976), 137-171.

\bibitem[H]{H} G. Hochschild, \emph{The structure of Lie groups}. Holden-Day, San Francisco, 1965.

\bibitem[LS]{LS} R. G. Larson and M. Sweedler, \emph{An associative orthogonal bilinear form for Hopf algebras}. Amer. J.  Math. \textbf{91} (1969), 75--94.

\bibitem[L]{L} B. I-Peng Lin, \emph{Semiperfect coalgebras}. J. Algebra {\bf 49} (1977), 357-373.

\bibitem[M]{M} S. Montgomery, \emph{Hopf algebras and their actions on rings}. CBMS Regional Conference Series in Mathematics {\bf 82},
Amer. Math. Soc., 1993.

\bibitem[R1]{Rad} D. E. Radford, \emph{On the coradical of a finite-dimensional Hopf algebra}. Proc. Amer. Math. Soc. \textbf{53} (1975), 9-15.

\bibitem[R2]{R} \bysame, {\it Finiteness conditions for a Hopf algebra with a nonzero integral}. J. Algebra {\bf 46} (1977), 189-195.

\bibitem[Sw1]{Sw-book} M. E. Sweedler, \emph{Hopf algebras}. Mathematics Lecture Note Series. Benjamin, 1969.

\bibitem[Sw2]{Sw} \bysame, \emph{Integrals for Hopf algebras}. Ann. of Math. (2) {\bf 89} (1969), 323-335.

\bibitem[T]{T} M. Takeuchi, \emph{Morita Theorems for Categories of Comodules.} J. Fac. Sci. Univ. Tokyo {\bf 24} (1977), 629-644.
\end{thebibliography}
\end{document}